\documentclass{scrartcl}

\usepackage{amsmath}
\usepackage{graphicx}
\usepackage{wrapfig}

\usepackage[standard,thmmarks,amsmath]{ntheorem}
\theoremstyle{plain}
\newtheorem{Assumption}{Assumption}

\usepackage{hyperref}

\usepackage[capitalise,noabbrev]{cleveref}
\crefname{Assumption}{Assumption}{Assumptions}

\usepackage{scalerel}
\usepackage{tikz}
\usetikzlibrary{svg.path}

\definecolor{orcidlogocol}{HTML}{A6CE39}
\tikzset{
  orcidlogo/.pic={
    \fill[orcidlogocol] svg{M256,128c0,70.7-57.3,128-128,128C57.3,256,0,198.7,0,128C0,57.3,57.3,0,128,0C198.7,0,256,57.3,256,128z};
    \fill[white] svg{M86.3,186.2H70.9V79.1h15.4v48.4V186.2z}
                 svg{M108.9,79.1h41.6c39.6,0,57,28.3,57,53.6c0,27.5-21.5,53.6-56.8,53.6h-41.8V79.1z M124.3,172.4h24.5c34.9,0,42.9-26.5,42.9-39.7c0-21.5-13.7-39.7-43.7-39.7h-23.7V172.4z}
                 svg{M88.7,56.8c0,5.5-4.5,10.1-10.1,10.1c-5.6,0-10.1-4.6-10.1-10.1c0-5.6,4.5-10.1,10.1-10.1C84.2,46.7,88.7,51.3,88.7,56.8z};
  }
}

\newcommand\orcidicon[1]{\href{https://orcid.org/#1}{\mbox{\scalerel*{
\begin{tikzpicture}[yscale=-1,transform shape]
\pic{orcidlogo};
\end{tikzpicture}
}{|}}}}

\begin{document}

\title{Machine learning of discrete field theories with guaranteed convergence and uncertainty quantification}
\author{Christian Offen ${\protect \orcidicon{0000-0002-5940-8057}}$\\
	{\small Paderborn University, Department of Mathematics}\\
	{\small Warburger Str. 100, 33098 Paderborn, Germany}\\
	{\small christian.offen@uni-paderborn.de}}

\def\d{\mathrm{d}}
\def\D{\mathrm{D}}
\def\p{\partial}
\def\R{\mathbb{R}}
\def\N{\mathbb{N}}
\def\Z{\mathbb{Z}}
\def\rf{{{\mathrm{ref}}}}
\def\DEL{{{\mathrm{DEL}}}}
\def\cntr{{{\mathrm{cntr}}}}

\maketitle              

\begin{abstract}
We introduce a method based on Gaussian process regression to identify discrete
variational principles from observed solutions of a field theory. The method is based on the data-based identification of a discrete Lagrangian density. It is a geometric machine learning technique in the sense that the variational structure of the true field theory is reflected in the data-driven model by design.
We provide a rigorous convergence statement of the method.
The proof circumvents challenges posed by the ambiguity of discrete Lagrangian densities in the inverse problem of variational calculus.
Moreover, our method can be used to quantify model uncertainty in the equations of motions and any linear observable of the discrete field theory.
This is illustrated on the example of the discrete wave equation and Schrödinger equation.
The article constitutes an extension of our previous article for the data-driven identification of (discrete) Lagrangians for variational dynamics from an ode setting to the setting of discrete pdes.
\end{abstract}

\section{Introduction}\label{sec:Introduction}

Partial differential equations that can be derived from variational principles (field theories) take a prominent role in physics, molecular biology, and engineering, for instance, as they describe wave phenomena, evolution of plasma dynamics, electro-magnetism, fluid dynamics, and quantum mechanical processes. Examples include the Korteweg--De Vries (KdV) equation, the shallow wave equation, and the Schrödinger equation.

In this context, data-driven technology has been used to obtain solutions to field theories \cite{Karniadakis2021} as well as to discover governing equations from observed solutions \cite{Qin2020,DLNNPDE}.
This article falls into the latter category. 
We propose to identify an action functional for a field theory by identifying a discrete Lagrangian density based on Gaussian process regression.
We do not assume any prior knowledge of the specific form of the Lagrangian. Central novelty of the article is the provision of a rigorous convergence theory for the proposed method. Moreover, Lagrangian descriptions of field theories are highly ambiguous \cite{HENNEAUX198245,Marmo1987,MARMO1989389}. This has practical implications for our machine learning framework that we discuss systematically and account for in the convergence theory. Furthermore, we provide systematic uncertainty quantification of linear observables of the data-driven system. 
The article can be seen as an extension of our article \cite{DLGPode} from the context of ordinary differential equations to a partial differential equations' context.

\paragraph{Continuous Lagrangian data-driven models}
A dynamical system is governed by a {\em variational principle} or a {\em least action principle}, if motions constitute critical points of an action functional. In case of an autonomous first-order field theory in the space of $\R^n$, the action functional is of the form
\begin{equation}\label{eq:SContinuous}
	S(u) = \int L(u(t),u_{t_1},\ldots,u_{t_n}) \d t_1 \ldots \d t_n,
\end{equation}
where $u$ is a scalar field and $u_{t_k} = \frac{\p u}{\p t_k}$ denote its derivatives ($k=1,\ldots,n$). The function $L$ is referred to as a Lagrangian density of the field theory.
In the special case of $n=1$ the function $L$ is called Lagrangian.
A function $u \colon \R^n \to \R$ is a solution of the field theory if for any bounded, open domain $\mathcal M \subset \R^n$ with Lebesgue measurable boundary $\p \mathcal{M}$ the action $S$ is stationary at $u|_{\mathcal{M}}$ for all variations $\delta u \colon \mathcal{M} \to \R$ that fix $\p \mathcal{M}$. This is (under regularity assumptions) equivalent to the condition that $u$ fulfils the Euler-Lagrange equations $\mathrm{EL}(u)=0$, where
\begin{equation}\label{eq:ELContinuous}
	\begin{split}
	\mathrm{EL}(u) &= \sum_{k=1}^{n} \frac{\p }{\p t_k} \left( \frac{\p L}{\p {u_{t_k}}} \right) - \frac{\p L}{\p u}
	= \sum_{k,l=1}^{n} \left( \frac{\p^2 L}{\p {u_{t_k}}\p {u_{t_l}}} u_{t_k,t_l} \right) - \frac{\p L}{\p u}.
	\end{split}
\end{equation}
Here, $u_{t_k,t_l} = \frac{\p^2 u}{\p t_k \p t_l}$.
Details may be found in \cite{gelfand2000calculus,RoubicekCalculusofVariations}, for instance.

In the data-driven context, $L$ is sought as a function of $u$ and its gradient such that \eqref{eq:ELContinuous} is fulfilled at observed data points \[\mathcal{D} = \left\{(u(t^{(j)}),(u_{t_k}(t^{(j)}))_{k=1}^n,(u_{t_k,t_l}(t^{(j)}))_{k,l=1}^n)\right\}_{j=1}^M\]
where $t^{(j)} \in \R^n$ for $1\le j \le M$. Once $L$ is known, $\mathrm{EL}(u)=0$ can be solved with a numerical method such as a variational integrator \cite{MarsdenWestVariationalIntegrators}.

In the special case of Lagrangian odes ($n=1$), several methods have been developed for the data-driven identification of a Lagrangian \cite{LNN,LagrangianShadowIntegrators,evangelisticdc2022,SymLNN,DLGPode}.

\paragraph{Discrete Lagrangian data-driven models}

Instead of learning continuous variational principles, in \cite{Qin2020} Qin proposes to learn discrete field theories by approximating discrete Lagrangian densities.
Here training data consists of solutions of the field theory observed over discrete spacetime locations, i.e.\ $\mathcal{D} = \{u(t^{(j)})\}_j$, where $t^{(j)}$ takes values in a temporal-spatial mesh.
This needs to be contrasted with the identification of a continuous model for a Lagrangian density which requires observations of derivatives of first and second order of solutions.

\begin{wrapfigure}{r}{0.28\textwidth} 
	\centering
	\includegraphics[width=\linewidth]{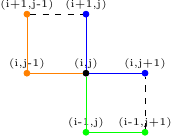}
	\caption{Stencil of \eqref{eq:7PtStencil}}\label{fig:7ptStencil}
\end{wrapfigure}
In case of a temporal-spatial mesh on $\R^2$ ($n=2$) with coordinates $(t,x)$, functions over the mesh (discrete fields) can be represented by the collection of their values $u^i_j$ at mesh points with indices $(i,j)$.
An example of a discrete Lagrangian density is the three-point Lagrangian $L_d(u^i_j,u^{i+1}_j,u^{i}_{j+1})$. A corresponding discrete action is defined as the (formal) sum
\begin{equation}\label{eq:SdExample}
S_\Delta(\{u^i_j\})
= \sum_{i,j} L_d(u^i_j,u^{i+1}_j,u^{i}_{j+1}) \Delta t \Delta x,
\end{equation}
where $(\Delta t, \Delta x)$ are the discretisation parameters of the mesh.
The application of a discrete variational principle yields the discrete Euler--Lagrange equation
\begin{equation}\label{eq:7PtStencil}
\frac{\p }{\p u^i_j}
\left(
L_d(u^i_j,u^{i+1}_j,u^{i}_{j+1})
+L_d(u^{i-1}_j,u^{i}_j,u^{i-1}_{j+1})
+L_d(u^{i}_{j-1},u^{i+1}_{j-1},u^{i}_{j})
\right)=0
\end{equation}
which corresponds to $\frac{\p }{\p u^k_l} S_\Delta(\{u^i_j\}) =0$ for all interior meshpoints $u^k_l$.
The discrete Euler--Lagrange equation constitutes a seven-point stencil which relates values of a discrete field over seven mesh-points as visualised in \cref{fig:7ptStencil}.
\begin{wrapfigure}{r}{0.28\textwidth} 
	\centering
	\includegraphics[width=\linewidth]{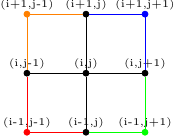}
	\caption{Stencil of \eqref{eq:9PtStencil}}\label{fig:9ptStencil}
\end{wrapfigure}
Another example is the four-point Lagrangian density $L_d(u^i_j,u^{i+1}_j,u^{i}_{j+1},u^{i+1}_{j+1})$. It yields the discrete Euler--Lagrange equation
\begin{equation}\label{eq:9PtStencil}
\begin{split}
\frac{\p }{\p u^i_j}
	\Big(
	&L_d(u^i_j,u^{i+1}_j,u^{i}_{j+1},u^{i+1}_{j+1})
	+L_d(u^{i-1}_j,u^{i}_j,u^{i-1}_{j+1},u^{i}_{j+1})\\
	+&L_d(u^{i}_{j-1},u^{i+1}_{j-1},u^{i}_{j},u^{i+1}_{j})
	+L_d(u^{i-1}_{j-1},u^{i}_{j-1},u^{i-1}_{j},u^{i}_{j})
	\Big)=0.
\end{split}
\end{equation}
This constitutes a nine-point stencil (\cref{fig:9ptStencil}). We refer to \cite{Marsden2001} for details on discrete variational calculus.
Three-point and four-point Lagrangian densities for discrete field theories over a temporal-spatial mesh will be the main examples in this article.


\paragraph{Ambiguity of Lagrangian densities and data-driven identification}
The data-driven identification of a continuous or a discrete Lagrangian density is an ill-defined inverse problem as many different Lagrangian densities can yield equations of motions with the same set of solutions \cite{HENNEAUX198245,Marmo1987,MARMO1989389}. In the context of neural network based learning methods, the author develops regularisation strategies that optimise numerical conditioning of the learnt theory in \cite{DLNNPDE,DLNNDensity}. The present article will provide normalisation conditions suitable for the context of learning Lagrangian densities with Gaussian processes.

The proposed method can be contrasted to operator learning (e.g. \cite{Lu2021}), which aims to learn an operator that can reconstruct solutions of pdes from boundary data.
Other approach aim to identify an explicit, analytic form of the governing pde via sparse regression based on a dictionary of possible expressions \cite{Schaeffer2017,Rudy2017}.
Furthermore, the method can be contrasted to methods based on model order reduction that seek a latent space on which the dynamics can be described by an ode that can be identified from data or derived from additional knowledge about a full order model. Structure-preserving approaches include \cite{Sharma2022,Sharma2022LagrangianROM,sharma2023symplectic, Blanchette2020,Blanchette2022}.

Our method exploits the temporal-spatial locality of partial differential equations by learning a comparatively low dimensional object, a discrete Lagrangian density, that can recover a discrete pde. The discrete pde can subsequently be solved for unseen boundary conditions by propagating data via the stencil defined by the discrete Euler--Lagrange equation. Moreover, our method preserves the temporal-spatial variational structure of the problem, which is related to multisymplectic conservation laws. In combination with variational symmetries variational structure is related to conservation laws via Noether's theorem \cite{Marsden1998}. Furthermore, preservation of local variational structure in machine learning approaches can be beneficial for the detection of travelling wave solutions \cite{DLNNPDE,DLNNDensity}.


\paragraph{Novelty of article}
The article constitutes an extension of \cite{DLGPode} from ordinary differential equations to partial differential equations.

\begin{enumerate}

\item We introduce a method based on Gaussian process regression to learn discrete Lagrangian densities from temporal-spatial data.

\item We formulate and prove a convergence theorem as the distance of data points tends to zero. 

\item The article systematically discusses the ambiguity of Lagrangian densities and normalisation strategies for kernel-based machine learning methods.

\item Moreover, a statistical framework for a quantification of model uncertainty of any linear observable of the learnt field theory is provided. 

\end{enumerate}

Our method may be contrasted to the aforementioned approaches in the literature for learning Lagrangian densities, which mainly focus on ordinary differential equations (potentially in combination with model order reduction) or on less geometric approaches such as operator learning. Moreover, ambiguity of sought Lagrangians is often circumvented in the literature by making a less general ansatz for the form of the Lagrangian (e.g.\ \cite{Aoshima2021}).
Additionally, convergence guarantees are typically not provided and model uncertainty quantification of linear observables is not discussed systematically, especially in the presence of model ambiguity (e.g.\ \cite{LNN}).

Interpreting a discrete Euler--Lagrange equation, such as \eqref{eq:7PtStencil} or \eqref{eq:9PtStencil}, as a partial differential equation for $L_d$, the methodology of our method can be viewed in the context of meshless collocation methods \cite{SchabackWendland2006} for linear partial differential equations. However, in contrast to approaches such as \cite{OwhadiLearningPDEGP}, solutions for $L_d$ of discrete Euler--Lagrange equations are highly ambiguous, which is a major technical difficulty that the article overcomes. Indeed, for this, we employ techniques presented by Owhadi and Scovel in \cite{OwhadiScovel2019OptimalRecoverySplines} to relate posterior means of Gaussian processes to constraint optimisation problems in reproducing kernel Hilbert spaces.


\paragraph{Outline}
In \cref{sec:background} the article continues the review of variational principles for discrete field theories and discusses the ambiguity of discrete Lagrangian densities for the description of field theories.
The analogous discussion for continuous theories is provided in \cref{app:ContinuousLagrangianDensities}, for comparison.
\Cref{sec:background}, furthermore, provides some observations that justify that our normalisation strategy does not restrict the generality of the ansatz. 
After a brief review of the notion of reproducing kernel Hilbert spaces and Gaussian processes, \cref{sec:MLSetting} introduces our method to identify Lagrangian densities from discrete temporal-spatial data with the possibility to quantify model uncertainty of linear observables.
\Cref{sec:Experiment} contains numerical experiments including the identification of a discrete Lagrangian density for the discrete wave equation and the Schrödinger equation.
\Cref{sec:ConvergenceAnalysis} contains a statement of convergence guarantees for the method.
The article concludes with a summary and concluding remarks in \cref{sec:Summary}.

As additional material, \cref{app:DiscreteLagrangiansFromPDETheorem} relates the convergence statement of \cref{sec:ConvergenceAnalysis} to the convergence analysis for temporal discrete Lagrangians presented in \cite{DLGPode}.

\section{Lagrangian dynamics}\label{sec:background}

In the following we continue our review of discrete Lagrangian dynamics for field theories. Here we focus our discussion on the three-point Lagrangian and four-point Lagrangian on a temporal-spatial mesh as introduced in \cref{sec:Introduction}.
The discussion can be expanded to other discrete Lagrangian densities that are obtained by applying quadrature formulas to exact discrete Lagrangian densities as defined in Section II.B of \cite{DLNNPDE}.
For details on discrete variational calculus we refer to \cite{Marsden2001}.
An analogous discussion for continuous theories can be found in \cref{app:ContinuousLagrangianDensities} for comparison.

\subsection{Discrete Euler--Lagrange equation}

As motivated in \cref{sec:Introduction} (see \eqref{eq:7PtStencil} and \eqref{eq:9PtStencil}), we define the discrete Euler--Lagrange operator for three-point Lagrangians by

\begin{equation}\label{eq:7PtDEL}
	\mathrm{DEL}(L_d)(\mathfrak{u})
	= \nabla_1 L_d(u,u^+,u_+)
	+\nabla_2 L_d(u^-,u,u^-_+)
	+\nabla_3 L_d(u_-,u^+_-,u)
\end{equation}
and for four-point Lagrangians by
\begin{equation}\label{eq:9PtDEL}
	\begin{split}
		\mathrm{DEL}(L_d)(\mathfrak{u})
		&=\nabla_1 L_d(u,u^+,u_+,u^+_+)
		+\nabla_2 L_d(u^-,u,u^-_+,u_+)\\
		&+\nabla_3 L_d(u_-,u^+_-,u,u^+)
		+\nabla_4 L_d(u^-_-,u_-,u^-,u).
	\end{split}
\end{equation}
Here $\mathfrak u = (u,u^+,u_+,u^-,u^-_+,u_-,u_-^+) \in (\R^d)^7$ for three-point Lagrangians and $\mathfrak u = (u,u^+,u_+,u_+^+,u^-,u^-_+,u_-,u_-^+,u_-^-) \in (\R^d)^9$ for four-point Lagrangians. The expression $\nabla_j L_d$ denotes the gradient of $L_d$ with respect to the $j$th input argument.
The discrete Euler--Lagrange equation is
\begin{equation}\label{eq:DELGeneral}
\mathrm{DEL}(L_d)(\mathfrak{u}) =0
\end{equation}
for all stencils $\mathfrak{u}$ that are contained in the mesh.


\subsection{Ambiguity of discrete Lagrangian densities}

The remainder of \cref{sec:background} prepares a justification that our normalisation strategy in the data-driven framework of \cref{sec:MLSetting} does not restrict the generality of the ansatz. 
Readers mostly interested in the data-driven aspects may skip directly to \cref{sec:MLSetting}.

Discrete Lagrangian densities can be ambiguous in two different ways: 
\begin{enumerate}
	
	\item Discrete Lagrangian densities $L_d$ and $\tilde L_d$ yield the same discrete Euler--Lagrange operator up to rescaling, i.e.
	\begin{equation}\label{eq:GaugeEqDefDEL}
		\mathrm{DEL}(L_d) = \rho \cdot  \mathrm{DEL}(\tilde{L}_d), \quad \rho \in \R\setminus \{0\}
	\end{equation}
	if
	\[
	L_d - \rho \tilde L_d -c = \mathrm{div}_t F,
	\]
	for $c \in \R$ and a discrete divergence $\mathrm{div}_t F$.
	For three-point Lagrangians $L_d(a,b,c)$ discrete divergence takes the form
	\[
	\mathrm{div}_t F(a,b,c)
	=F_1(a)-F_1(b) + F_2(a)-F_2(c)+F_3(b)-F_3(c)
	\]
	with continuously differentiable function $F=(F_1,F_2,F_3)$. These terms are telescopic in the action $S_d$ \eqref{eq:SdExample}.
	For four-point Lagrangians $L_d(a,b,c,d)$ a discrete divergence takes the form
	\begin{equation*}
		\begin{split}
			\mathrm{div}_t F(a,b,c,d)
			&=
			F_1(a)-F_1(b) 
			+ F_2(a)-F_2(c)
			+ F_3(a)-F_3(d)\\
			&+ F_4(b)-F_4(c) 
			+ F_5(b)-F_5(d)
			+ F_6(c)-F_6(d)	
		\end{split}
	\end{equation*}
	with continuously continuously differentiable function $F=(F_1,F_2,F_3,F_4,F_5,F_6)$.
	Two Lagrangians $L$ and $\tilde L$ with \eqref{eq:GaugeEqDefDEL} are called {\em (gauge-) equivalent}.
	
	\item Even when two discrete Lagrangian densities $L_d$ and $\tilde L_d$ are not equivalent, they can yield the same set of solutions, i.e.
	\[
	\mathrm{DEL}(L_d)(\mathfrak{u}) =0 \iff \mathrm{DEL}(\tilde L_d)(\mathfrak{u}) =0,
	\]
	where $\mathfrak{u}$ denotes stencil data. 
	In that case $\tilde{L}_d$ is called an {\em alternative Lagrangian density} to $L_d$.

\end{enumerate}

\subsection{Normalisation of discrete Lagrangian densities}\label{sec:Normalisation}

In the data-driven framework for learning discrete Lagrangian densities that we will introduce in a subsequent section, we will employ strategies to prevent degenerate solutions of discrete Lagrangians that fail to discriminate non-motions. An extreme example are Null-Lagrangians, for which $\mathrm{DEL}(L_d)\equiv 0$. 

To a four-point discrete Lagrangian density $L_d$
define $\mathrm{Mm}^+(L_d)$ by
\begin{equation}\label{eq:MmDef}
	\begin{split}
	\mathrm{Mm}^+(L_d)(\mathfrak{u}) &= \nabla_2 L_d(u^-,u,u^-_+,u_+)  + \nabla_4 L_d(u^-_-,u_-,u^-,u)\\
	\mathrm{Mm}^-(L_d)(\mathfrak{u}) &= -\nabla_1 L_d(u,u^+,u_+,u^+_+) -\nabla_3 L_d(u_-,u^+_-,u,u^+).
	\end{split}
\end{equation}
for elements of nine-point stencils $\mathfrak{u} = (u,u^+,u_+,u^-,u^-_+,u_-,u_-^+) \in (\R^d)^9$.
Considering discrete three-point Lagrangian densities as four-point Lagrangian densities with trivial dependence on the fourth component, the definitions extend to three-point Lagrangian densities.

\begin{remark}[Interpretation]
	The expression $\mathrm{Mm}^+(L_d)$ and $\mathrm{Mm}^-(L_d)$ can be viewed as a component of a conjugate momentum: to a discrete four-point Lagrangian density $L_d$ consider the Lagrangian
	\begin{equation}\label{eq:LdDX}
	L_{d, \Delta x}(U^i,U^{i+1}) = \sum_j L_d(u^i_j,u^{i+1}_j,u^{i}_{j+1},u^{i+1}_{j+1}),
	\end{equation}
	where $U^i = (u^i_j)_j$. The discrete Lagrangian describes the temporal motion of the system via the discrete Euler--Lagrangian equation for ordinary differential equations
	\[
	0 = \nabla_{U^i} \left( L_{d, \Delta x}(U^{i-1},U^{i}) + L_{d, \Delta x}(U^i,U^{i+1}) \right).
	\]
	Associated discrete conjugate momenta are defined by $P_+^i = \nabla_{U^i}  L_{d, \Delta x}(U^{i-1},U^{i})$ and $P_-^i = - \nabla_{U^i} L_{d, \Delta x}(U^{i},U^{i+1})$. Indeed, the $j$th component is given as
	\begin{align*}
		(P_+^i)_j &= \nabla_2 L_d(u^{i-1}_j,u^{i}_j,u^{i-1}_{j+1},u^{i}_{j+1})  + \nabla_4 L_d(u^{i-1}_{j-1},u^{i}_{j-1},u^{i-1}_{j},u^{i}_{j})\\
		(P_-^i)_j &= -\nabla_1 L_d(u^i_j,u^{i+1}_j,u^{i}_{j+1},u^{i+1}_{j+1}) -\nabla_3 L_d(u^i_{j-1},u^{i+1}_{j-1},u^{i}_{j},u^{i+1}_{j}).
	\end{align*}

\end{remark}

\begin{proposition}\label{prop:normaliseLd}
	Let $\mathfrak{u}_b = (u,u^+,u_+,u^-,u^-_+,u_-,u_-^+) \in (\R^d)^9$ be a stencil and $\mathring{L_d}$ be a discrete four-point Lagrangian.
	To any $c_b \in \R$, $p_b \in \R^d$ there exists an equivalent discrete Lagrangian $L_d$ with
	\begin{equation*}\label{eq:NormaliseLd}
		L_d(u,u^+,u_+,u^+_+) = c_b,
		\qquad
		\mathrm{Mm}^-(L_d)(\mathfrak{u}_b) = p_b.
	\end{equation*}
\end{proposition}

\begin{proof}
	Let $\mathring c_b = \mathring{L_d}(u,u^+,u_+,u^+_+)$, $\mathring{p_b}=\mathrm{Mm}^-( \mathring{L_d})(\mathfrak{u}_b)$, $F=(F_1,0,0,0,0,0)$ with
	$F_1(u) = (\mathring{p_b}-p_b)^\top u$, and $c = c_b - \mathring{c_b} + (\mathring{p}_b - p_b)(u^+-u)$.
	The Lagrangian
	\[
	L_d = \mathring{L_d} + \mathrm{div}_t F + c
	\]
	is equivalent to $\mathring{L_d}$ and fulfils the conditions \eqref{eq:NormaliseLd}.
\end{proof}

\begin{remark}
	An analogous statement holds with $\mathrm{Mm}^+$ replacing $\mathrm{Mm}^-$ in \cref{prop:normaliseLd}.
\end{remark}

\section{Data-driven framework}\label{sec:MLSetting}

Consider an open and bounded domain $\Omega \subset (\R^d)^q$ ($q \in \N$) on which an (a priori unkown) discrete Lagrangian density $L_d^\rf \colon \Omega \to \R$ defines true motions via $\mathrm{DEL}(L_d^\rf) (\mathfrak{u}) =0$.
Based on the observation of stencil data $\mathfrak{u}$ for which $\mathrm{DEL}(L_d^\rf) (\mathfrak{u}) =0$, our goal is to approximate a Lagrangian density $L_d \colon \Omega \to \R$ such that 
$\mathrm{DEL}(L_d^\rf) (\mathfrak{u}) =0 \iff \mathrm{DEL}(L_d) (\mathfrak{u}) =0$ for all stencils $\mathfrak u \in \Omega$.

For this, we will model $L_d$ as a Gaussian Process over $\Omega$ and compute the posterior process of $L_d$ conditioned on $\mathrm{DEL}(L_d) (\mathfrak{u}) =0$ for all observed stencil data $\mathfrak{u}$. To prevent learning a degenerate $L_d$ (such as Null-Lagrangians) we will, furthermore, condition on additional normalisation conditions as motivated by \cref{prop:normaliseLd}.

For this, we briefly introduce our set-up making use of a reproducing kernel Hilbert-space (RKHS) setting. We will then proceed to the computation of the posterior process.
We refer the reader to \cite{ChristmannSteinwart2008RKHS,OwhadiScovel2019} for an introduction and further information.

\subsection{Reproducing kernel Hilbert spaces -- Set-up}

Consider a continuously differentiable function $K \colon \Omega \times \Omega \to \R$. Assume that $K$ is a positive-definite kernel, i.e.\ $K(x,y)=K(y,x)$ for $x,y \in \Omega$ (symmetry) and the matrix $(K(x^{(i)},x^{(j)}) )_{i,j=1}^M$ is positive definite for all finite subsets $\{x^{(i)}\}_{i=1}^M \subset \Omega$.

To define the reproducing kernel Hilbert space $U$ to $K$, 
consider the set $\{K(x^{(j)},\cdot )\; | \; x^{(j)} \in \Omega \}$ and its linear span $\mathring{U}$. An inner product on $\mathring{U}$ is provided by the linear extension of $\langle \cdot, \cdot \rangle$ defined by
\[
\langle K(x, \cdot), K(y, \cdot) \rangle = K(x,y).
\]
The reproducing kernel Hilbert space $U$ is now obtained as the completion of $\mathring{U}$ with respect to $\langle \cdot, \cdot \rangle$.

An identification of the dual space $U^\ast$ and $U$ is given by the bijective, linear, symmetric map $\mathcal{K} \colon U^\ast \to U$ defined by $\mathcal{K}(\Phi)(x) = \Phi(K(x, \cdot))$.

Let $\xi \in \mathcal{N}(0, \mathcal{K})$ be the {\em canonical Gaussian process} on $U$, i.e.\ $\xi \colon \mathcal{A} \to U$ is a random variable (where $\mathcal{A}$ is a probability space) with $\mathbb{E}(\xi) =0$ (zero mean) and covariance operator $\mathcal{K}$ as introduced above. This means that for $\Phi = (\Phi_1,\ldots,\Phi_m) \in (U^\ast)^m$ ($m \in \N$) the random variable $\Phi(\xi)$ on $\R^m$ is multivariate-normally distributed with covariance matrix given by $\kappa = (\Phi_i (\mathcal{K}(\Phi_j)))_{i,j =1}^m$, i.e.\ $\Phi(\xi) \in \mathcal{N}(0,\kappa)$.

\subsection{Computation of the posterior Gaussian Process}\label{sec:PosteriorComputation}

In case we are seeking a three-point discrete Lagrangian density, let us assume we are given stencil data
\[\mathfrak{u}^{(k)} = (u^{(k)},
{u^{(k)}}^+,
{u^{(k)}}_+,
{u^{(k)}}^-,
{u^{(k)}}^-_+,
{u^{(k)}}_-
{u^{(k)}}^+_-
) \in (\R^d)^7\]
for $k=1,\ldots,M$
such that the sub-triples are elements of $\Omega$, i.e.\
$({u^{(k)}},{u^{(k)}}^+,{u^{(k)}}_+)$,
$({u^{(k)}}^-,{u^{(k)}},{u^{(k)}}^-_+)$,
$({u^{(k)}}_-,{u^{(k)}}^+_-,{u^{(k)}}) \in \Omega$.

In case a four-point discrete Lagrangian density is sought, we consider stencil data of the form
\[\mathfrak{u}^{(k)} = ({u^{(k)}},{u^{(k)}}^+,{u^{(k)}}_+,{u^{(k)}}_+^+,{u^{(k)}}^-,{u^{(k)}}^-_+,{u^{(k)}}_-,{u^{(k)}}_-^+,{u^{(k)}}_-^-) \in (\R^d)^9\]
for $k=1,\ldots,M$
such that the sub-quadruples are elements of $\Omega$, i.e.\
$({u^{(k)}},{u^{(k)}}^+,{u^{(k)}}_+,{u^{(k)}}^+_+)$,
$({u^{(k)}}^-,{u^{(k)}},{u^{(k)}}^-_+,{u^{(k)}}_+)$,
$({u^{(k)}}_-,{u^{(k)}}^+_-,{u^{(k)}},{u^{(k)}}^+)$,
$({u^{(k)}}^-_-,{u^{(k)}}_-,{u^{(k)}}^-,{u^{(k)}}) \in \Omega$.

Moreover, consider an additional stencil $\mathfrak{u}_b \in (\R^d)^7$ or $\mathfrak{u}_b\in (\R^d)^9$, respectively, for which the sub-triples (sub-quadruples) are elements of $\Omega$ and values $p_b \in \R^d$, $c_b \in \R$ to which we define $y^M_b := (0,\ldots,0,p_b,c_b) \in (\R^d)^M \times \R^d \times \R$.

Moreover, let $\Phi_b^M \colon U \to (\R^d)^M \times \R^d \times \R$ be defined as
\begin{equation}\label{eq:PhibM}
	\Phi_b^M = (\mathrm{DEL}_{\mathfrak{u}^{(1)}},\ldots,\mathrm{DEL}_{\mathfrak{u}^{(M)}},\mathrm{Mm}^-_{\mathfrak{u}_b},\mathrm{ev}_{\mathfrak{u}_b}),
\end{equation}
where for $L_d \in U$
\begin{align*}
	\mathrm{DEL}_{\mathfrak{u}^{(k)}}(L_d) &:= \mathrm{DEL}(L_d)(\mathfrak{u}^{(k)}), \quad 1\le k \le M\\
	\mathrm{Mm}^-_{\mathfrak{u}_b}(L_d)
	&:=\mathrm{Mm}^-(L_d)(\mathfrak{u}_b)\\
	\mathrm{ev}_{\mathfrak{u}_b}(L_d)
	:=L_d(u_b,u_b^+,{u_b}_+,{u_b}_+^+) &\text{ or } \mathrm{ev}_{\mathfrak{u}_b}(L_d) := L_d(u_b,u_b^+,{u_b}_+)
\end{align*}
with $\mathrm{Mm}^-$ as defined in \eqref{eq:MmDef}.


The conditional process $\xi^M := \xi | \{\Phi^M_b(\xi)=y^M_b\}$ is again a Gaussian process $\xi^M \in \mathcal N(L_d^M,\mathcal{K}_{\Phi^M_b})$ since $\Phi^M_b$ is linear.
Thus, $\xi^M$ is fully defined by the conditional mean $L_d^M \in U$ and conditional covariance operator $\mathcal{K}_{\Phi^M_b} \colon U^\ast \to U$. These can be computed using general theory \cite[Cor.~17.12]{OwhadiScovel2019} as follows.

Define the symmetric, matrix $\Theta \in \R^{(M+1)d \times (M+1)d}$ as
\begin{equation}\label{eq:ThetaMatLd}
	\Theta = \begin{pmatrix}
		(\mathrm{DEL}^1_{\mathfrak{u}^{(j)}}\mathrm{DEL}^2_{\mathfrak{u}^{(i)}}K)_{ij}
		&(\mathrm{DEL}^1_{\mathfrak{u}^{(j)}}\mathrm{Mm^-}^2_{\mathfrak{u}_b}K)_{j}
		&(\mathrm{DEL}^1_{\mathfrak{u}^{(j)}} \mathrm{ev}^2_{\mathfrak u_b}K)_j\\
		(\mathrm{Mm^-}^1_{\mathfrak{u}_b}\mathrm{DEL}^2_{\mathfrak{u}^{(i)}}K)_{i}
		&\mathrm{Mm^-}^1_{\mathfrak{u}_b}\mathrm{Mm^-}^2_{\mathfrak{u}_b}K
		&\mathrm{Mm^-}^1_{\mathfrak{u}_b} \mathrm{ev}^2_{\mathfrak u_b}K\\
		(\mathrm{ev}^1_{\mathfrak u_b}\mathrm{DEL}^2_{\mathfrak{u}^{(i)}}K)_{i}
		&\mathrm{ev}^1_{\mathfrak u_b}\mathrm{Mm^-}^2_{\mathfrak{u}_b}K
		&\mathrm{ev}^1_{\mathfrak u_b}\mathrm{ev}^2_{\mathfrak u_b}K.
	\end{pmatrix}
\end{equation}
The upper index $1,2$ of the operators $\mathrm{DEL}$, $\mathrm{Mm}^-$, $\mathrm{ev}$ denote on which input element of $K$ they act, i.e.\
\[
\mathrm{DEL}^1_{\mathfrak{u}^{(j)}}\mathrm{DEL}^2_{\mathfrak{u}^{(i)}}K
= \mathrm{DEL}_{\mathfrak{u}^{(j)}} (a\mapsto  \mathrm{DEL}_{\mathfrak{u}^{(i)}}(b \mapsto K(a,b))),
\]
for instance. Moreover, we have used the convention that operators are applied component-wise to vectors of functions.


The Gaussian process $\xi^M \in  \mathcal{N}(L_d,\mathcal{K}_{\Phi_b^M})$ has conditional mean 
\begin{equation}\label{eq:LdPosterior}
	L_d = {y_b^M}^\top \theta^{-1} \mathcal{K}\Phi_b^M.
\end{equation}

The conditional covariance operator $\mathcal{K}_{\Phi_b^M} \colon U^\ast \to U$ is defined for any $\psi,\phi \in U^\ast$ by
\begin{align}\label{eq:CovPosteriorLd}
	&\psi\mathcal{K}_{\Phi_b^M}\phi
	= \psi\mathcal{K}\phi
	- (\psi\mathcal{K}{\Phi_b^M}^\top) \theta^{-1} (\Phi_b^M\mathcal{K}\phi).
\end{align}
Here
\begin{align*}
	\psi\mathcal{K}_{\Phi_b^M}\phi	&=\psi^1 \phi^2 K\\
	\psi\mathcal{K}{\Phi_b^M}^\top &= \begin{pmatrix}
		\psi^1\mathrm{DEL}^2_{\mathfrak{u}^{(2)}}K,
		& \ldots &
		\psi^1\mathrm{DEL}^2_{\mathfrak{u}^{(M)}}K,
		&
		\psi^1\mathrm{Mm^-}^2_{\mathfrak{u}_b}K,
		& \psi^1 K(\cdot ,u_b)
	\end{pmatrix}\\
	\Phi_b^M\mathcal{K}\phi &=
	\begin{pmatrix}
		\mathrm{DEL}^1_{\mathfrak{u}^{(2)}}\phi^2K &
		\ldots &
		\mathrm{DEL}^1_{\mathfrak{u}^{(M)}}\phi^2K &
		\mathrm{Mm^-}^1_{\mathfrak{u}_b}\phi^2K &
		\phi^2 K(u_b,\cdot )
	\end{pmatrix}^\top.
\end{align*}


\begin{Remark}[Normalisation]
We would like that the discrete Euler--Lagrange equation \eqref{eq:DELGeneral} to a data-driven discrete Lagrangian density uniquely defines a solution over the temporal-spatial mesh for suitable boundary data. This is referred to as non-degeneracy of the Lagrangian density. (The exact notion  depends on the type of boundary value problem.)
The normalisation condition $(\mathrm{Mm}^-_{\mathfrak{u}_b}(\xi),\mathrm{ev}_{\mathfrak{u}_b}(\xi)) = (p_b,c_b)$
is compatible with the true motions defined by $L^\rf_d$ as the normalisation condition is fully covered by gauge freedom as proved in \cref{prop:normaliseLd}.
Moreover, as by \cite[Thm 12.5]{OwhadiScovel2019OptimalRecoverySplines}, the conditional mean $L_d$ can be characterised as the minimiser of the convex optimisation problem
\begin{equation}\label{eq:ConvexOpt}
	L_d = \mathrm{argmin}_{\Phi^M_b(\tilde L_d)=y_b^M} \|\tilde L_d\|_U,
\end{equation}
where $U$ is the reproducing kernel Hilbert space to the kernel $K$. 
Even though the normalisation does not guarantee any non-degeneracy of the conditional mean $L_d$, we will see in the numerical examples that the condition $p_b \not =0$, $c_b \not = 0$ causes well-behaved discrete Lagrangian densities for the prediction of new solutions to the field theory from unseen training data. Intuitively, discrete Lagrangians $L_d$ with artificial degeneracies do not happen to constitute minimisers of \eqref{eq:ConvexOpt} in practise unless the degeneracies are forced by the true underlying field theory.
\end{Remark}

\begin{Remark}[Computational aspects]
The evaluation of the posterior mean $L_d$ or covariance operator $\mathcal{K}_{\Phi_b^M}$ involves solving large linear systems of equations, see \eqref{eq:LdPosterior},\eqref{eq:CovPosteriorLd}.
While we will employ standard libraries in the numerical experiments, we refer to the literature on computational methods for Gaussian process regression for various approaches such as \cite{QuinRasmussen2005,Schaefer2021,ChenOwhadiSChaefer2024SparseCholesky}.
\end{Remark}

\section{Numerical experiments}\label{sec:Experiment}

In this section, we apply the method of \cref{sec:MLSetting} to observational data of the discrete wave equation and a discrete Schrödinger equation. Similar experiments have been performed in \cite{DLNNPDE} for neural network based methods.

\subsection{Discrete wave equation}\label{sec:WaveExperiment}

\subsubsection{True model and training data generation}

\begin{figure}
	\centering
	\includegraphics[width=0.32\linewidth]{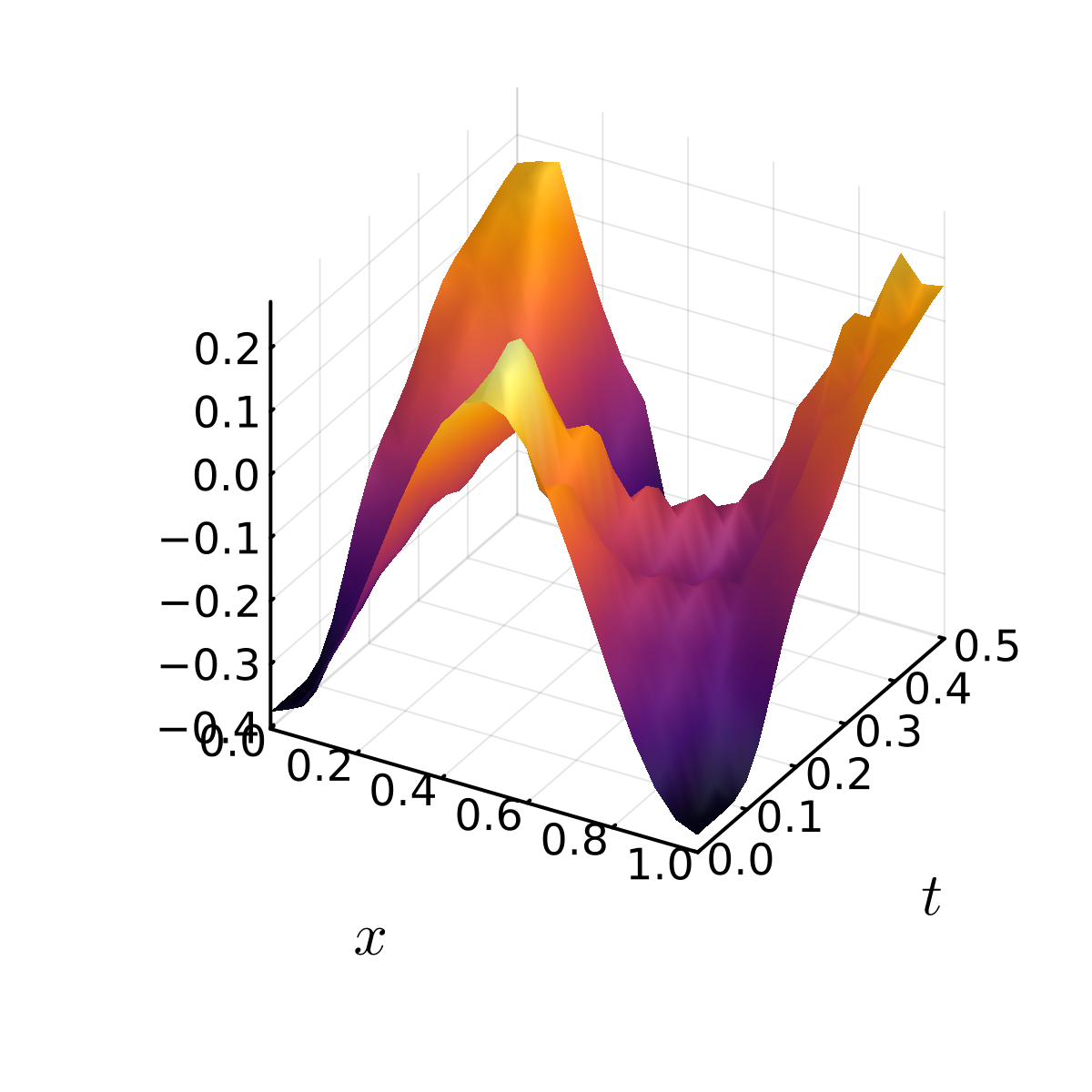}
	\includegraphics[width=0.32\linewidth]{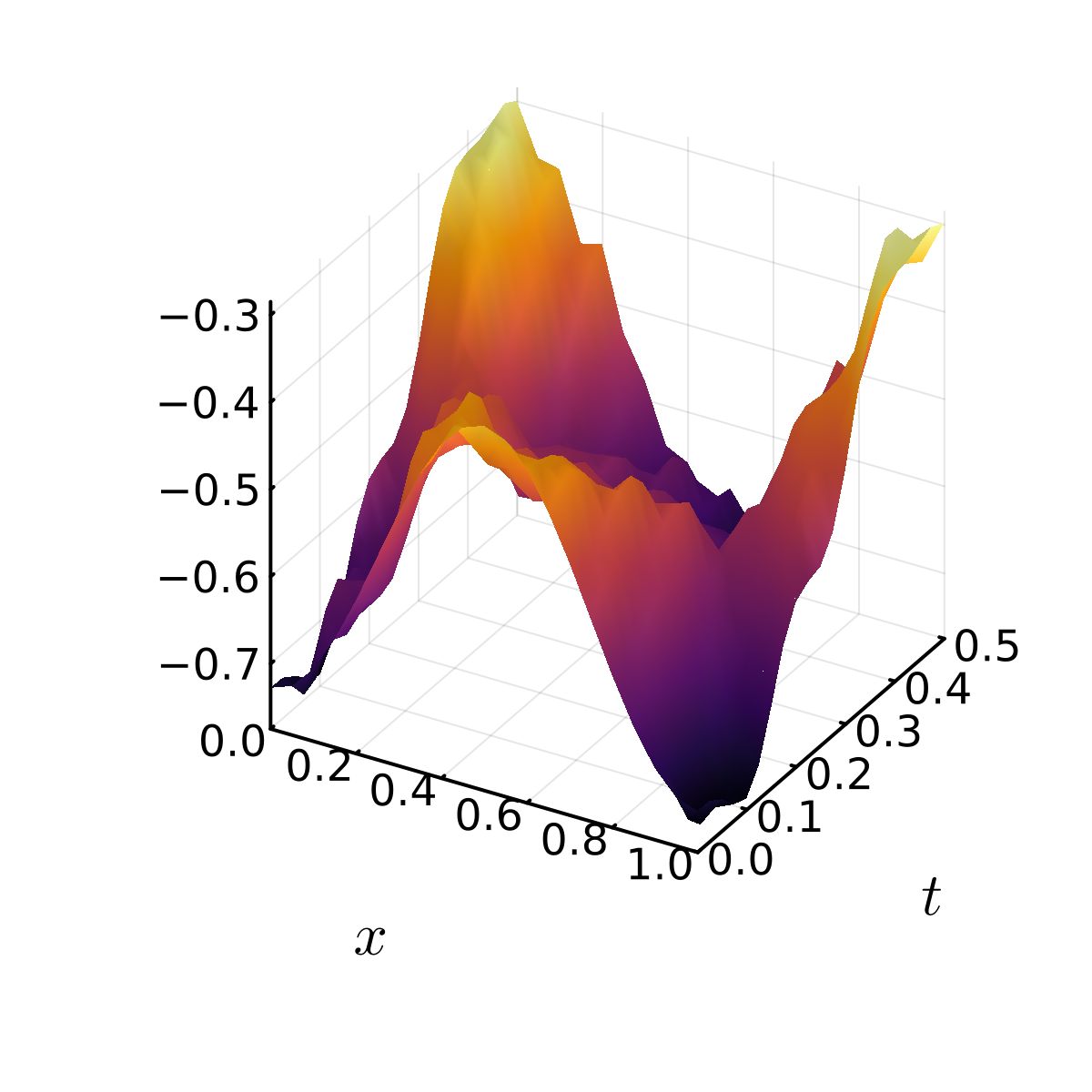}
	\caption{Complete training data set for wave equation experiment. (Only these two samples.)}\label{fig:SamplesWaveEQ}
\end{figure}

On the temporal-spatial domain $\Omega = [0,0.5] \times [0,1]$ with coordinates $(t,x)$ and periodic boundary conditions in space, consider a uniform mesh with discretisation parameters $\Delta t = 1/40$, $\Delta x = 1/20$ and the following discrete Lagrangian
\[
L^\rf_d(u,u^+,u_+) = \frac 12 \left(\frac{u^+-u}{\Delta t}\right)^2 - \frac 12 \left(\frac{u_+-u}{\Delta x}\right)^2 - V(u), \qquad \text{with }V(u) = \frac 12 u^2.
\]
The discrete Euler--Lagrange equations
\begin{equation}\label{eq:DiscreteWave}
	0=\mathrm{DEL}(L_d^\rf)(\mathfrak u)
	= \frac{u^--2u+u^+}{\Delta t^2} - \frac{u_--2u+u_+}{\Delta x^2} + \nabla V(u)
\end{equation}
is the discrete wave equation and can be viewed as a discrete version of the continuous wave equation $\frac {\p^2 u}{\p t^2}-\frac {\p^2 u}{\p x^2} + \nabla V(u) =0$.
To obtain training data, we sample $N_0=2$ initial conditions consisting of discrete spatial data $U^0 = (u^0_j)_j$ and discrete conjugate momentum data $P^0_- = (p^0_j)_j$ at $t=0$. The initial data are then propagated in time by forming the temporal Lagrangian density $L_{d,\Delta x}$ defined in \eqref{eq:LdDX}. Then standard methods in variational integration \cite{MarsdenWestVariationalIntegrators} are used:
in the first step, $P_-^0 = - \nabla_{U^0} L_{d, \Delta x}(U^{0},U^{1})$ is solved for $U^1$. Then $U^{i+1}$ is computed from $U^{i}$ and $U^{i-1}$ by solving the discrete Euler--Lagrange equation $\mathrm{DEL}(L^\rf_{d,\Delta x})(U^{i-1},U^{i},U^{i+1}) =0$ numerically for $U^{i+1}$ ($i\ge 1$). The discrete Euler--Lagrange operator $\mathrm{DEL}$ for temporal Lagrangians is given as
\begin{equation}\label{eq:DELLdx}
\mathrm{DEL}(L_{d,\Delta x})(U^{i-1},U^{i},U^{i+1})
= \frac{\p }{\p U^i}
\left(
L_{d,\Delta x}(U^{i-1},U^{i})
+ L_{d,\Delta x}(U^{i},U^{i+1})
\right).
\end{equation}
The two samples are plotted in \cref{fig:SamplesWaveEQ}.
These two solutions contain $M=760$ seven-point stencils $\mathfrak U = (\mathfrak{u}^{(k)})_{k=1}^M$ ($N_0$ times number of interior mesh points).

\subsubsection{Predictions from unseen initial data (Extrapolation)}

With the base stencil $\mathfrak{u}_b = 0 \in \R^7$, normalisation constants $c_b=1$, $p_b=1$, and kernel $K(x,y) = \exp(-1/2 \| x-y\|^2)$, the stencil data $\mathfrak U$ defines a posterior Gaussian process $\xi^M \in \mathcal{N}(L_d^M,\mathcal{K}_{\Phi_b^M})$ (\cref{sec:PosteriorComputation}).

\Cref{fig:PredictWaveEQ} shows that the mean $L_d^M$ of the posterior process can successfully be used to predict solutions of the discrete wave equation from the unseen initial data $u(0,x)=-\cos(2 \pi x)$, $u(\Delta t,x)=-\cos(2 \pi x)$ via forward propagation with $\mathrm{DEL}(L_{d,\Delta x}^M)=0$. The discrete $l_2$-norm of the error is smaller than $0.0222$.

To quantify model uncertainty along the predicted solution, we consider at each interior mesh point $(i \Delta t,j \Delta x)$ the stencil $\mathfrak{u}^{(i,j)}$ centred at $(i \Delta t,j \Delta x)$. The standard deviation $\sigma$ of the Gaussian random variable $\mathrm{DEL}_{\mathfrak{u}^{(i,j)}}(\xi^M)$ can be computed using the conditional covariance operator $\mathcal K_{\Phi^M_b}$ as follows
\begin{equation}\label{eq:StndDEL}
\sigma(\mathrm{DEL}_{\mathfrak{u}^{(i,j)}}(\xi^M))
= \sqrt{\psi \mathcal K_{\Phi^M_b} \phi}, \text{ with } \psi = \phi = \mathrm{DEL}_{\mathfrak{u}^{(i,j)}}.
\end{equation}
Here $\psi \mathcal K_{\Phi^M_b} \phi$ is computed as explained in \eqref{eq:CovPosteriorLd}.
The plot to the right of \cref{fig:PredictWaveEQ} displays the values of $\sigma(\mathrm{DEL}_{\mathfrak{u}^{(i,j)}}(\xi^M))$ at all interior mesh points $(i,j)$.
This shows areas of very low and slightly increased model uncertainty and can be interpreted as an estimation of the local error.

\subsubsection{Prediction of travelling waves (Extrapolation)}

\begin{figure}
	\centering
	\includegraphics[width=0.32\linewidth]{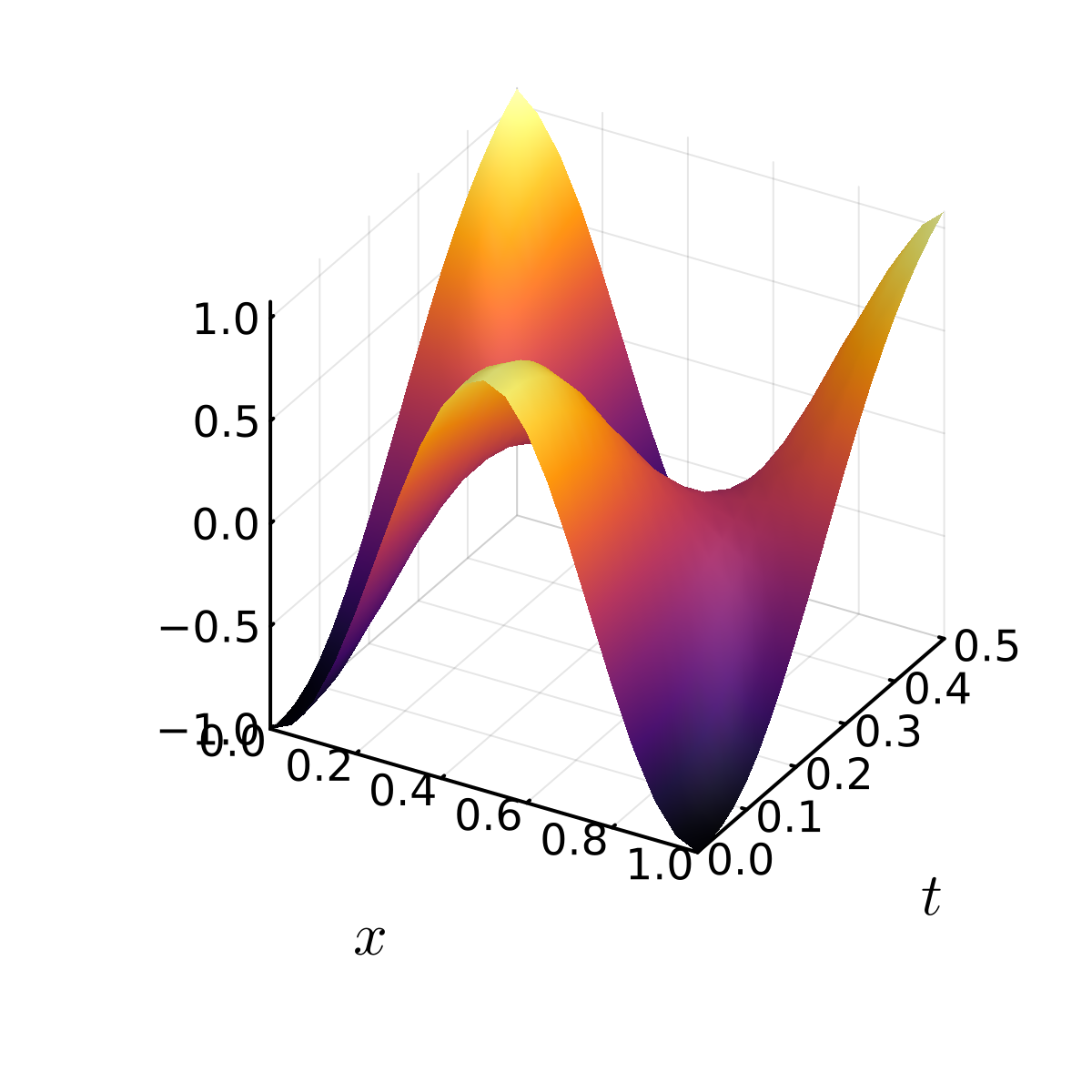}
	\includegraphics[width=0.32\linewidth]{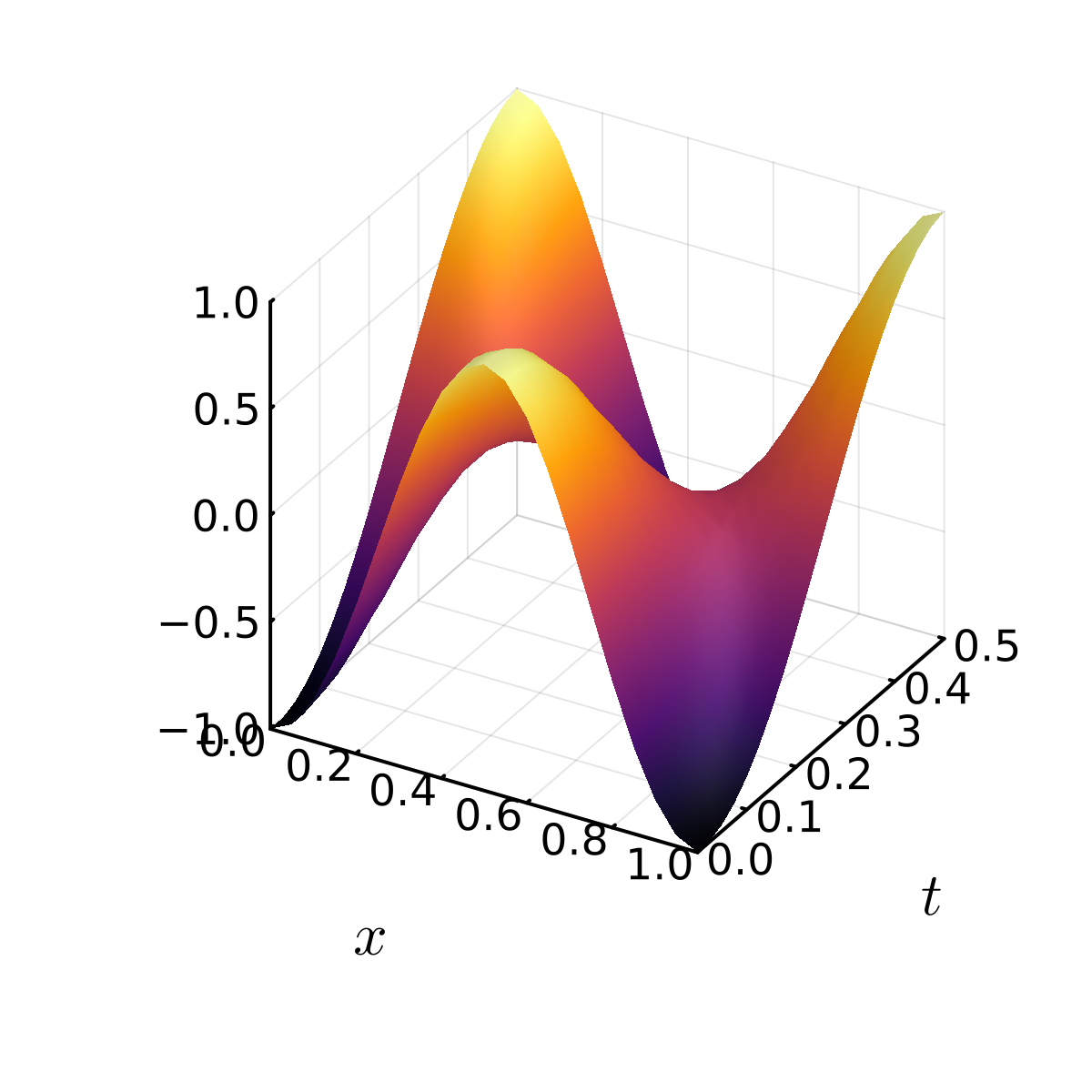}
	\includegraphics[width=0.3\linewidth]{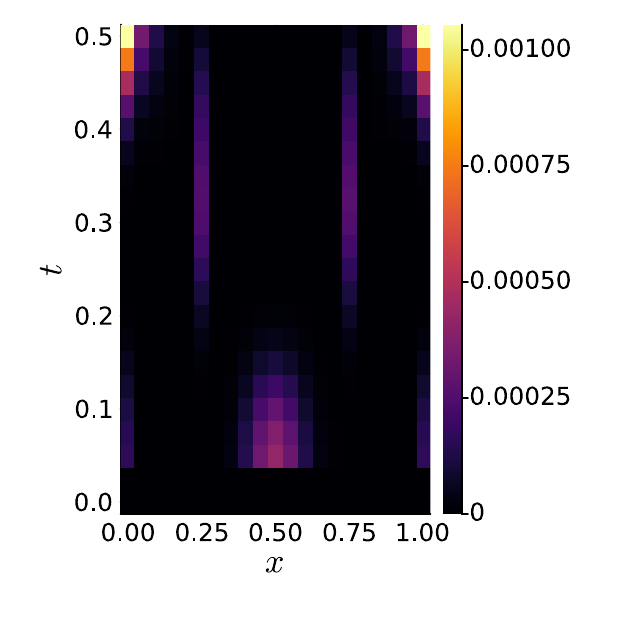}
	\caption{Wave equation experiment. Predicted solution from initial data at $t=0$ with $L_d^M$ (left) matches reference (centre) to high accuracy. Right: Standard deviation of the Gaussian variable $\mathrm{DEL}(\xi^M)$ along predicted solution.}\label{fig:PredictWaveEQ}
\end{figure}

\begin{figure}
	\centering
	\includegraphics[width=0.32\linewidth]{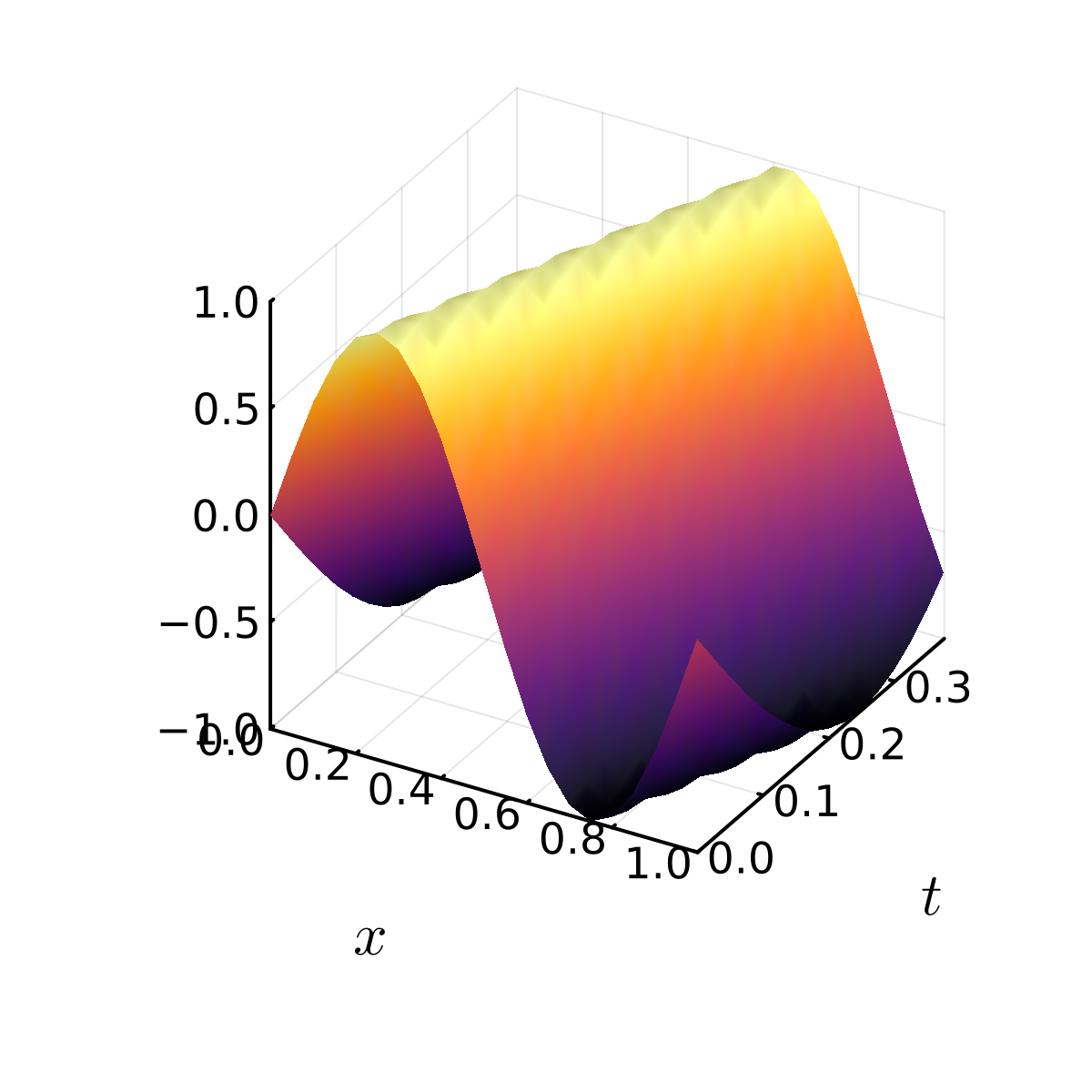}
	\includegraphics[width=0.32\linewidth]{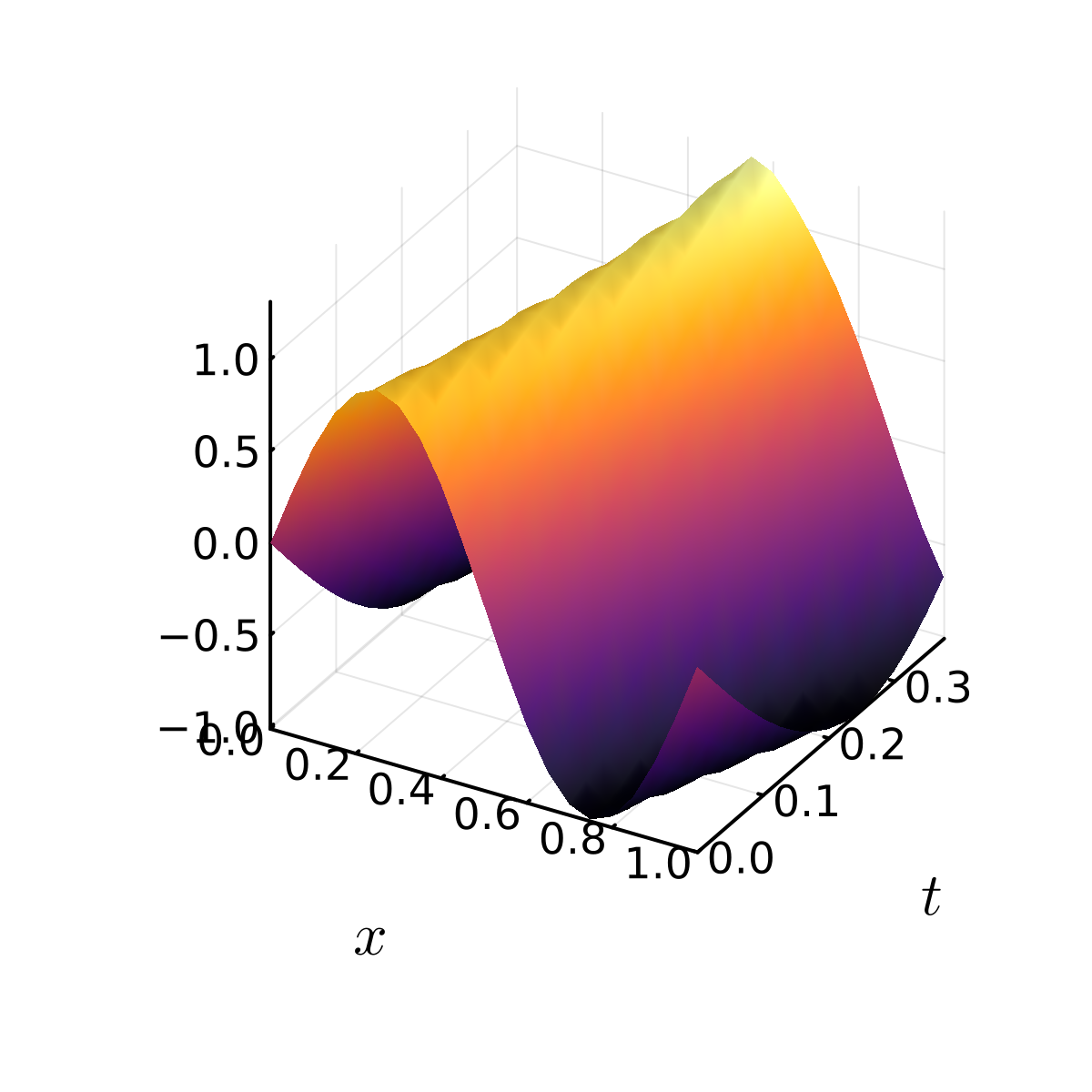}
	\includegraphics[width=0.3\linewidth]{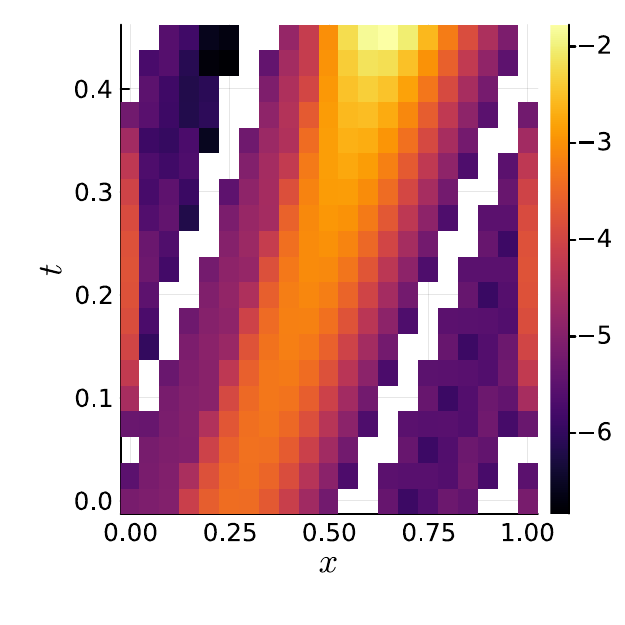}
	\caption{Wave equation experiment. Predicted solution from travelling wave initial data at $t=0$ with $L_d^M$ (left) matches reference (centre) to high accuracy but with an upward trend of the ridge. Right: Standard deviation of the Gaussian variable $\mathrm{DEL}(\xi^M)$ along predicted solution. For better visualisation, $\log_{10}$ is applied before plotting a heat map.}\label{fig:PredictWaveEQTW}
\end{figure}
The functional equation
\begin{align*}
	&\phantom{+}\frac{u(t-\Delta t,x)-2u(t,x)+u(t+\Delta t,x)}{\Delta t^2}\\
	&- \frac{u(t,x-\Delta x)-2u(t,x)+u(t,x+\Delta x)}{\Delta x^2}
	+ u(t,x)
	=0
\end{align*}
to a continuous interpolation of \eqref{eq:DiscreteWave}
admits travelling wave solutions $u_{\mathrm{TW}}(t,x) = f_{\mathrm{TW}}(s) = a_1 \sin(\kappa s) + a_2 \cos(\kappa s)$, $s=x-ct$, $\kappa = 2\pi k$ ($k \in \mathbb{Z}$, $a_1,a_2 \in \R$), where the wave speed $c$ is a real solution to 
\[
\cos(\kappa c \Delta t) = 1 - \frac{\Delta t^2}{2} + \frac{\Delta t^2}{\Delta x^2}(\cos(k \Delta x)-1)
\]
(see \cite[Example 11]{DLNNPDE}).
\Cref{fig:PredictWaveEQTW} shows a repetition of the experiment of \cref{fig:PredictWaveEQ} for the initial data of a travelling wave $u_{\mathrm{TW}}(0,x)$, $u_{\mathrm{TW}}(\Delta t,x)$. An upward trend of the ridge can be observed. The discrete $l_2$-norm of the error is smaller than 0.067 (time interval $[0,0.375]$).
In the heatmap of the standard deviation (right plot), for better visualisation, the logarithm (basis 10) is applied to the data before plotting. The upward trend of the ridge is reflected in a higher model uncertainty at later times.

The experiments demonstrate that the learned model can extrapolate to unseen initial data and that enough information is contained in the stencil data of the two observed fields (\cref{fig:SamplesWaveEQ}) to predict travelling waves to high accuracy.
This is remarkable as the training data-set does not contain travelling waves. Indeed, the global behaviour of the two fields in the training data set \eqref{fig:SamplesWaveEQ}) is quite different from the global behaviour of travelling waves. However, as the considered technique for identification of the dynamical system and the method to predict solutions relies on local dependencies, extrapolation to solutions with different global behaviour is possible.

\subsection{Discrete Schrödinger equation}\label{sec:ExperimentSE}

\subsubsection{True model and training data generation}

Let $\Omega = [0,0.14] \times [0,1]$ be a temporal-spatial domain with coordinates $(t,x)$ and periodic boundary conditions in $x$. On the domain consider a uniform mesh that is periodic in space with mesh parameters $(\Delta t, \Delta x)$. To a smooth potential $V \colon \R \to \R$ consider the discrete Lagrangian
\[
L_d^\rf(u,u^+,u_+,u^+_+) =
u_{\cntr}^\top J u_{\Delta t}
- \|u_{\Delta x} \|^2
- V(\|u_{\cntr} \|^2)
\]
with
\begin{equation*}
	\begin{split}
		J&=\begin{pmatrix}	0 & -1\\1 & 0\end{pmatrix},\\ 
		u_{\cntr} &= \frac 1 4(u + u^+ + u_+ + u^+_+),
	\end{split}
	\qquad 
	\begin{split}
		u_{\Delta t} &= \frac 1 {2 \Delta t} ((u^+_+ - u_+) + (u^+ - u) ),\\
		u_{\Delta x} &= \frac 1 {2 \Delta x} ((u^+_+ - u^+) + (u_+ - u) ).
	\end{split}
\end{equation*}
The discrete Lagrangian yields the following discrete Euler--Lagrange equation
\begin{equation}
	\begin{split}
		0&=\mathrm{DEL}(L_d^\rf)(\mathfrak{u})
		= 2Ju_{\Delta t,\cntr} + 2 u_{\Delta x^2} - 2 \left(V'(\|u_\cntr\|^2)u_\cntr\right)_\cntr\\ \label{eq:DELSE}
		\iff Ju_{\Delta t,\cntr} &= u_{\Delta x^2} -  \left(V'(\|u_\cntr\|^2)u_\cntr\right)_\cntr
	\end{split}
\end{equation}
where we have used the notation
\begin{align*}
	u_{\Delta t,\cntr} &= \frac{1}{8 \Delta t}((u^+_+-u^+_-)+2(u_+-u_-) + (u^-_+-u^-_-))\\
	u_{\Delta x^2} &= \frac{1}{4 \Delta x^2}
	(
	(u^+_+-2u^++u^+_-)
	+2(u_+-2u+u_-)
	+(u_+^--2u^-+u_-^-))\\
	\left(g(u,u^+,-u_+,u^+_+)\right)_\cntr
	&=\frac 14 (\mathrm{evaluate}+\mathrm{shift}^-+\mathrm{shift}_-+\mathrm{shift}^-_-)\left(g,(u,u^+,u_+,u^+_+)\right).
\end{align*}
Above, the operator $\mathrm{shift}$ evaluates $g$ on shifted input arguments. The shift occurs in the temporal or spatial direction as indicated by its indices.
The discrete Euler--Lagrange equation \eqref{eq:DELSE} constitutes a discrete version of the continuous Schrödinger equation
\begin{equation}\label{eq:SEComplex}
	\hbar \mathrm{i} \frac{\p }{\p t} \Psi = \left(- \frac{\p^2}{\p x^2} + V'(|\Psi|^2)\right)\Psi,
\end{equation}
which in real coordinates $u=(\mathrm{Re}(\Psi),\mathrm{Im}(\Psi))^\top$ reads
\begin{equation}\label{eq:SEReal}
	\hbar J \frac{\p }{\p t} u = \left(- \frac{\p^2}{\p x^2} + V'(\|u\|^2)\right)u.
	\end{equation}
(The Plank constant is set to $\hbar=1$.) Equation \eqref{eq:SEReal} may be interpreted as a continuous analogue to \eqref{eq:DELSE}.
For the purpose of the numerical experiments, the discrete equation \eqref{eq:DELSE} is regarded as the true field theory. We may refer to the first component of $u$ as real part, and to the second component as imaginary part.

In the experiments, $V(r) = r$, $\Delta t=7/400$, $\Delta x=1/10$ ($9\times 10$ mesh points).
To obtain training data, we sample $N_0=30$ initial spatial data at time $t=0$ and integrate forward in time using variational integration techniques as explained in \cref{sec:WaveExperiment}. However, the momentum $P^0$ is not sampled but computed from $U^0$ exploiting the special structure of $L_d^\rf$ 
(refer to \cite[II B.3]{DLNNPDE} for details).

\begin{figure}
	\centering
	\includegraphics[width=0.24\linewidth]{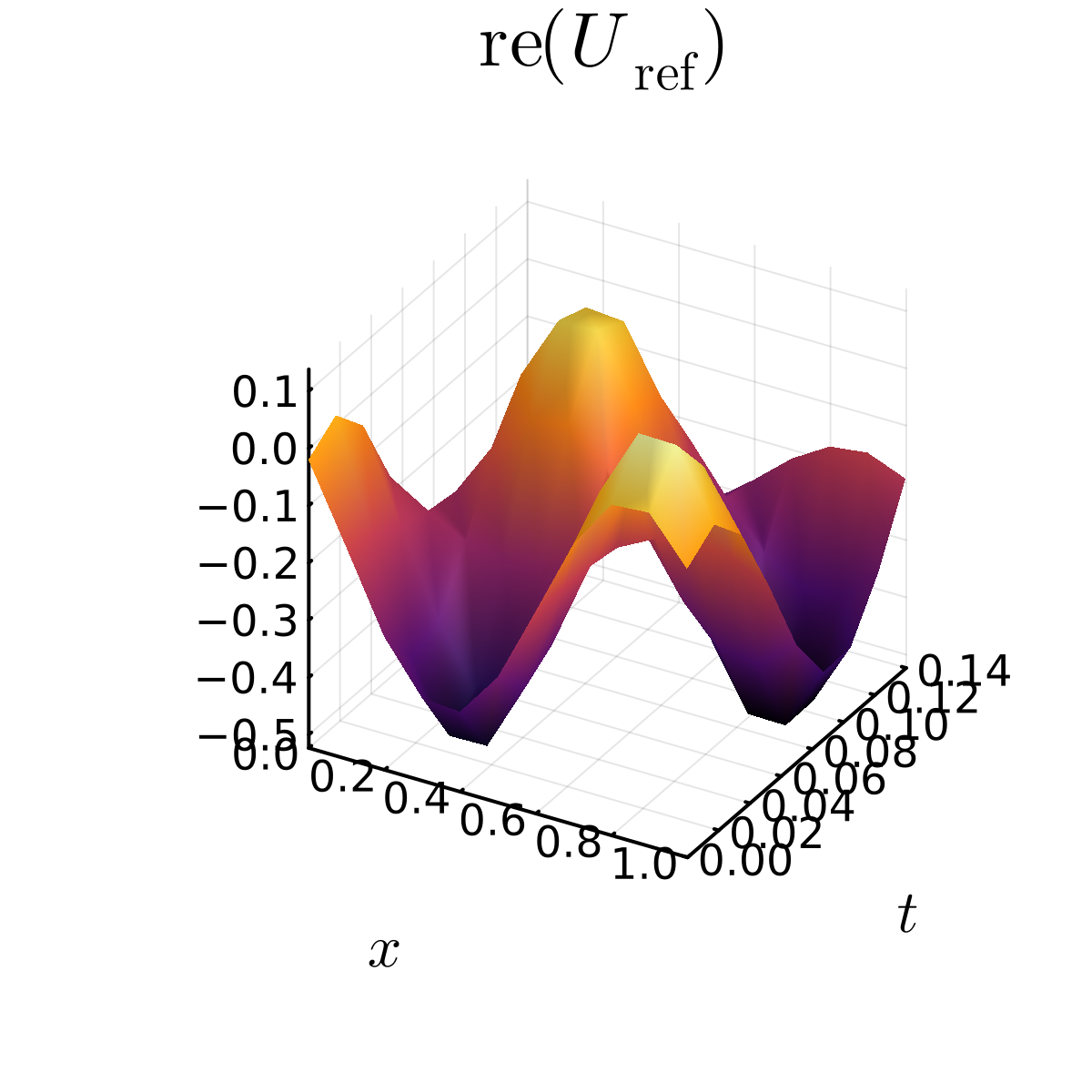}
	\includegraphics[width=0.24\linewidth]{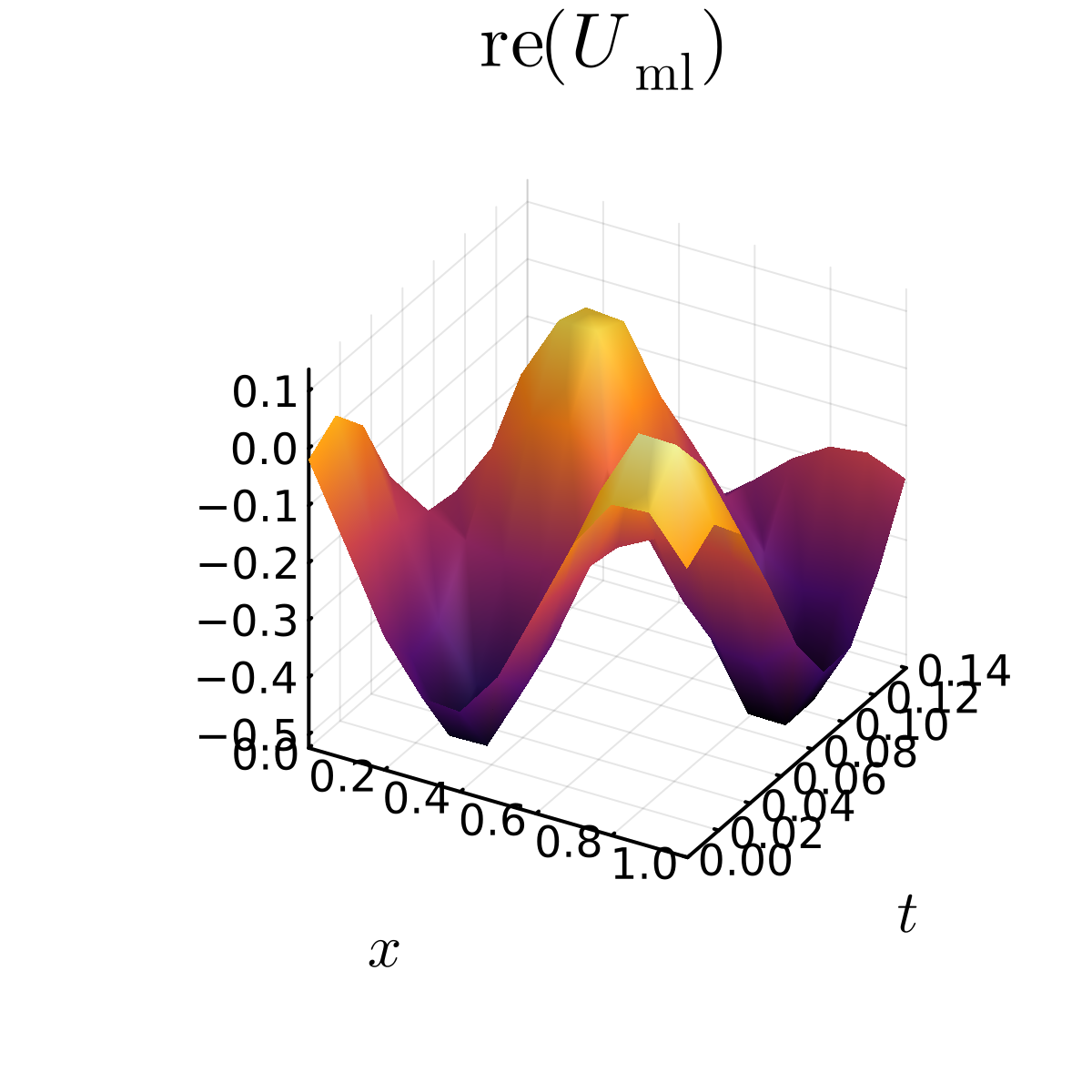}
	\includegraphics[width=0.24\linewidth]{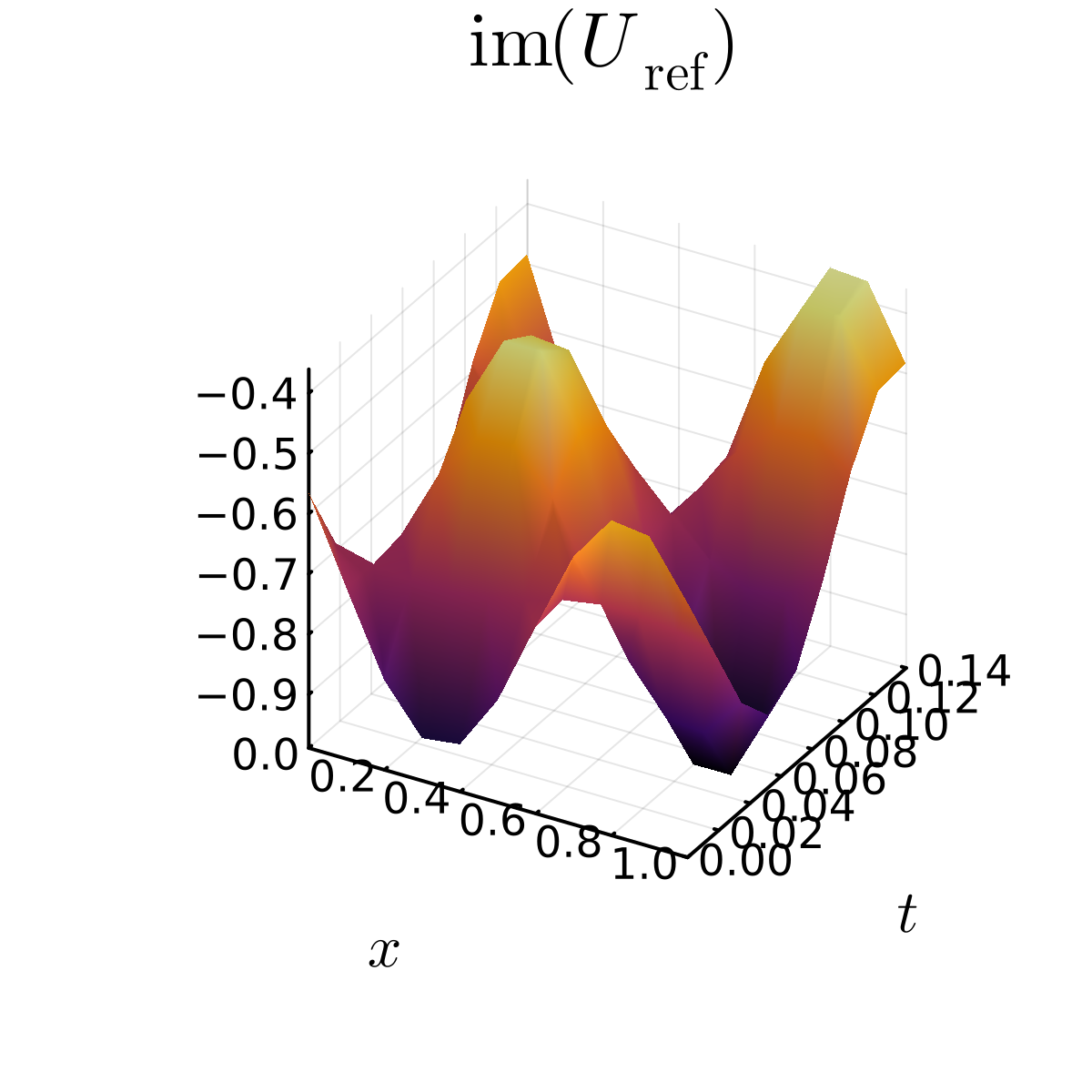}
	\includegraphics[width=0.24\linewidth]{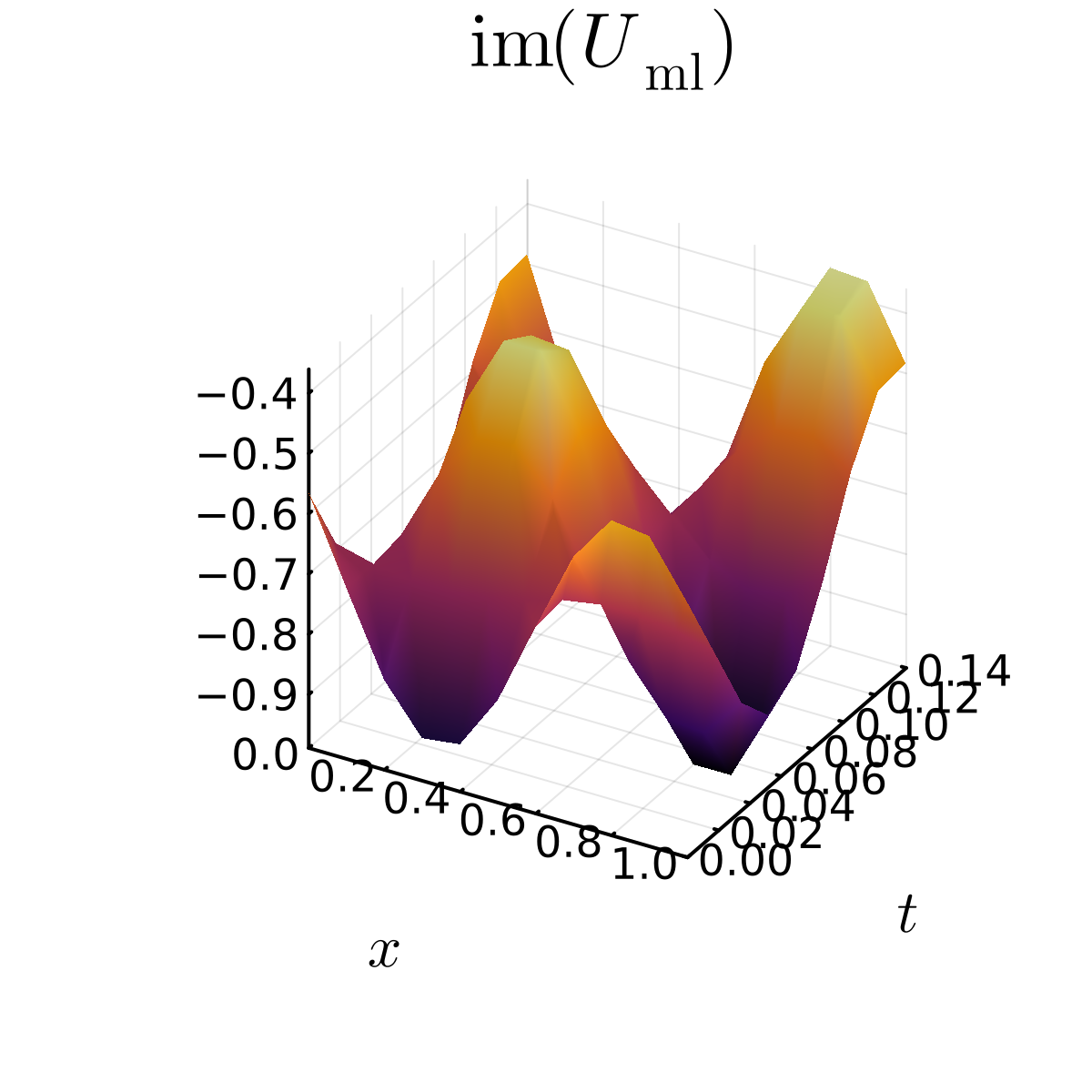}
	\caption{Schrödinger equation experiment.
		Recovery of training data from initial values at $t=0,\Delta t$ by the learned model.
		From left to right: real part of training data, real part of predicted field,
		imaginary part of training data, imaginary part of predicted field.
	}\label{fig:TrainingRecoverSE}
\end{figure}

The first and third plot of \cref{fig:TrainingRecoverSE} display real and imaginary part of an instance of a training sample. These solutions constitute training data $\mathfrak U =(\mathfrak u^{(k)})_{k=1}^M$ consisting of $M=2100$ nine-point stencils (one stencil per interior mesh point per sample).
With the base stencil $\mathfrak{u}_b = 0 \in (\R^2)^9$, normalisation constants $c_b=1$, $p_b=\begin{pmatrix}1,&1\end{pmatrix}^\top$, and kernel $K(x,y) = \exp(-1/2 \| x-y\|^2)$, the stencil data $\mathfrak U$ defines a posterior Gaussian process $\xi^M \in \mathcal{N}(L_d^M,\mathcal{K}_{\Phi_b^M})$ (see \cref{sec:PosteriorComputation}).

\subsubsection{Recovery of training data from initial values}

Neglecting effects of finite-precision arithmetic and assuming non-degeneracy, by construction, training data can be reproduced exactly from initial data
at $t=0,\Delta t$ by forming the temporal Lagrangian $L_{d,\Delta x}^M$ and forward propagation in time using $\mathrm{DEL}(L_{d,\Delta x}^M)=0$. Moreover, the model uncertainty $\sigma(\mathrm{DEL}_{\mathfrak u}(\xi^M))$ at stencils $\mathfrak u$ in the training data $\mathfrak U$ is zero in exact arithmetic.
This is confirmed in the numerical experiment in \cref{fig:TrainingRecoverSE}. The discrete $l_2$ norm of reference and predicted field is smaller than $1.1e-9$. The maximal value of each component of $\sigma(\mathrm{DEL}_{\mathfrak u}(\xi^M))$ is smaller than $1.45e-7$ for $\mathfrak{u}\in \mathfrak{U}$.

\subsubsection{Extrapolation}

\begin{figure}
	\centering
	\includegraphics[width=0.32\linewidth]{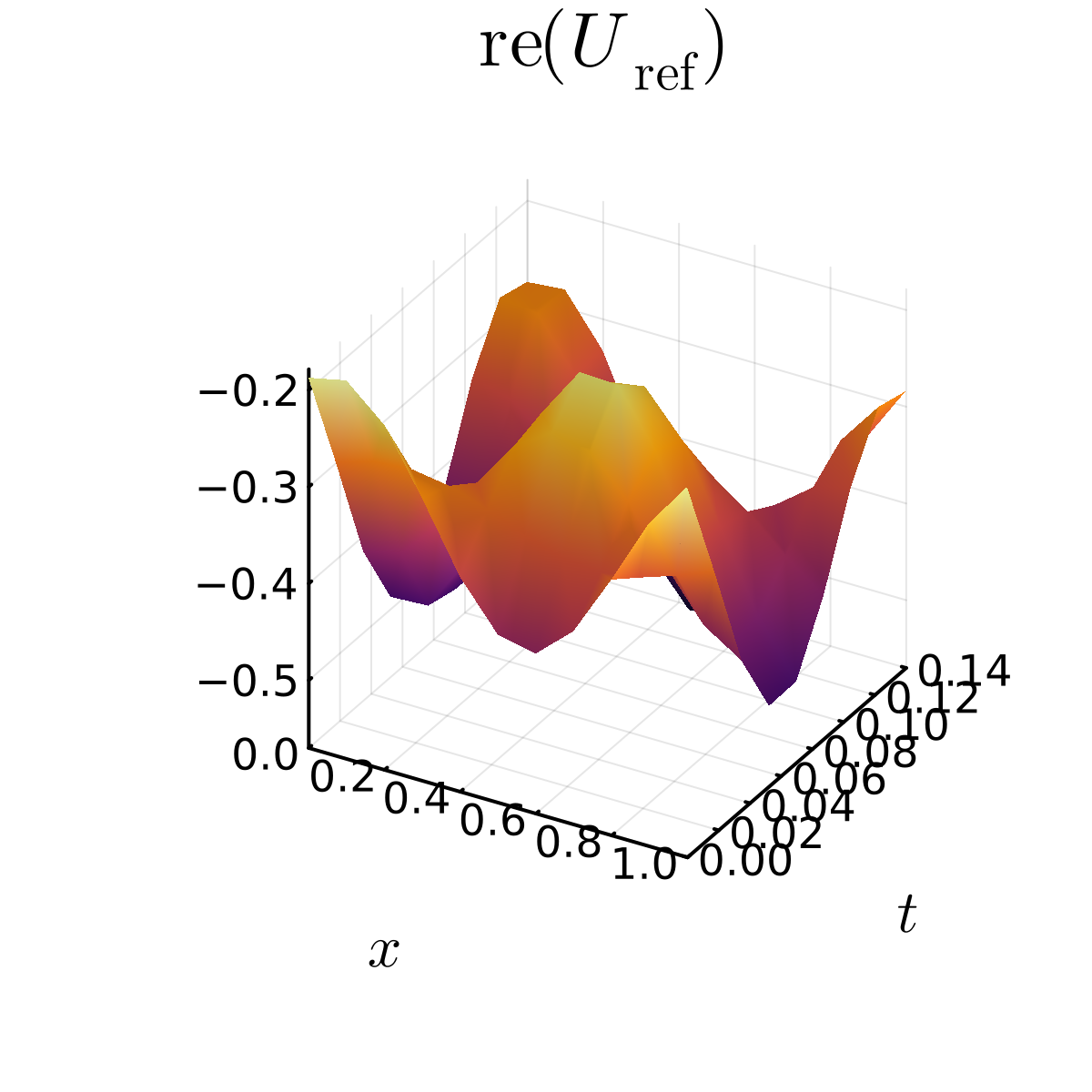}
	\includegraphics[width=0.32\linewidth]{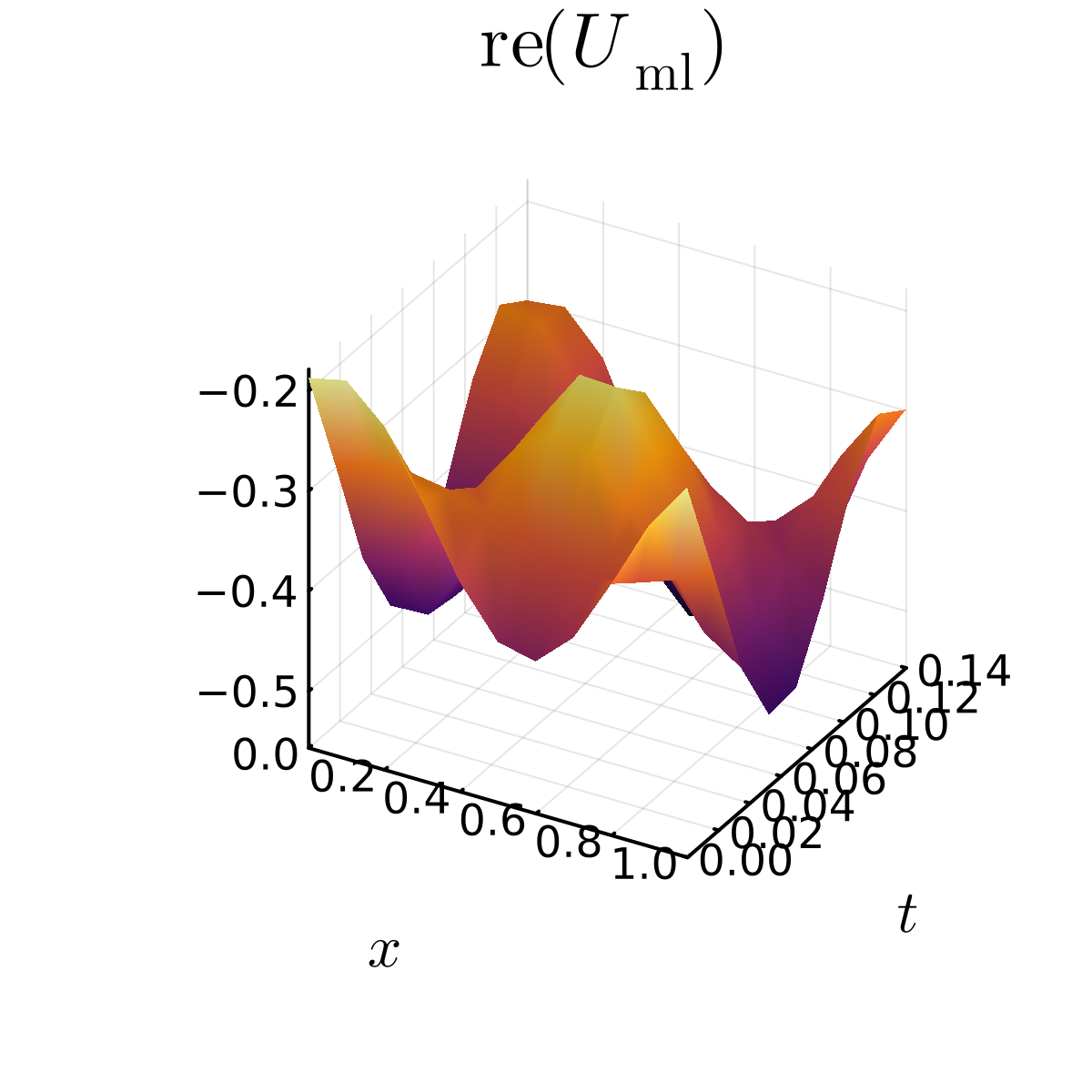}
	\includegraphics[width=0.32\linewidth]{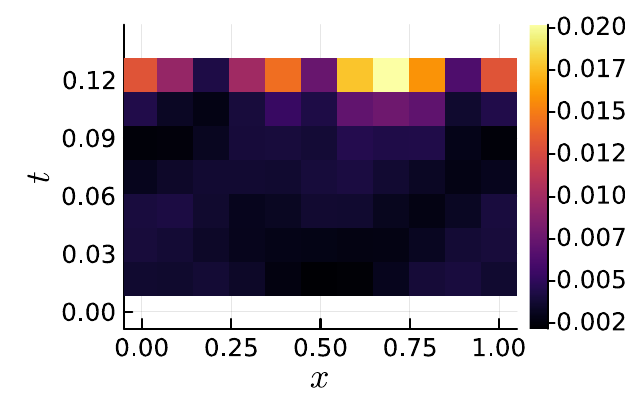}\\
	\includegraphics[width=0.32\linewidth]{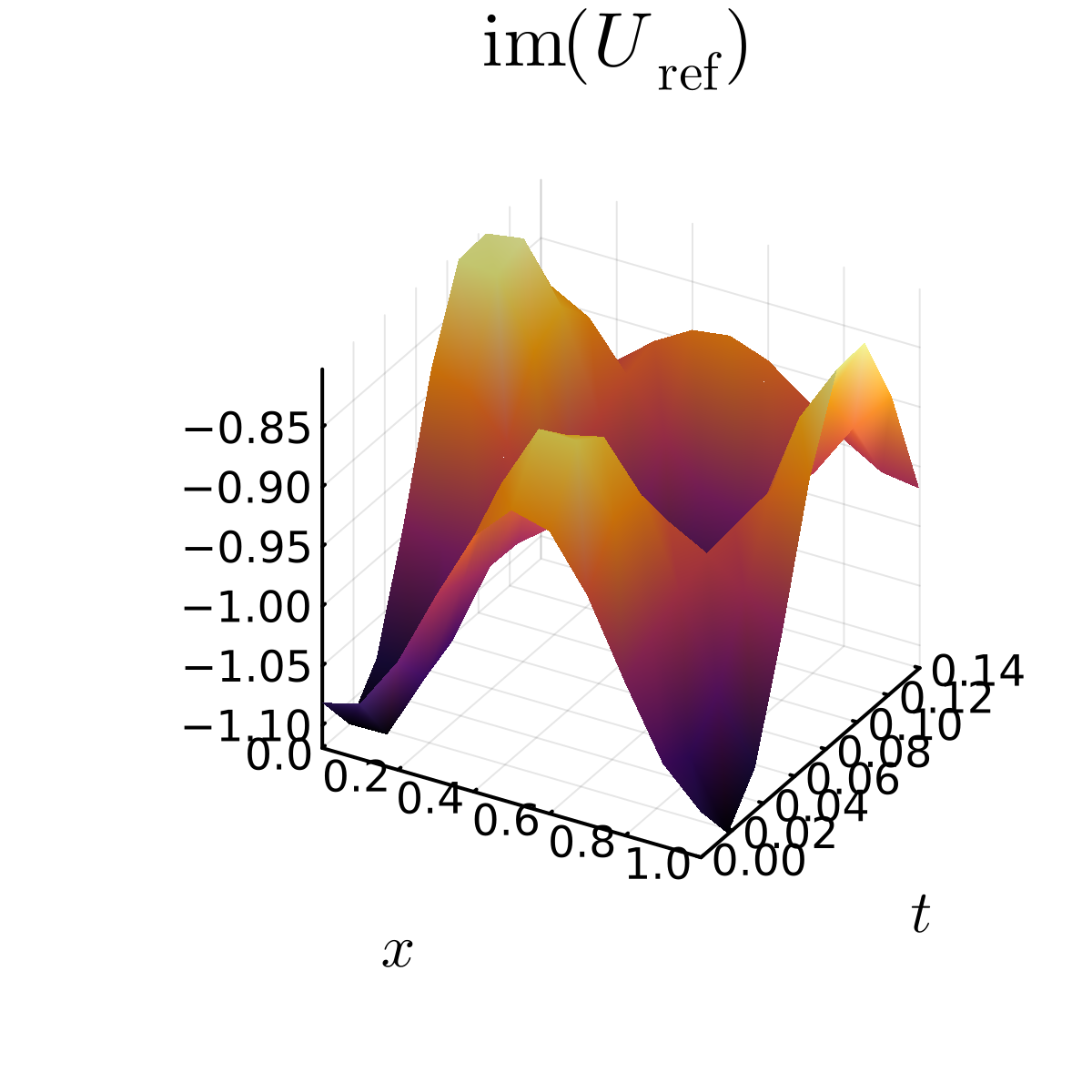}
	\includegraphics[width=0.32\linewidth]{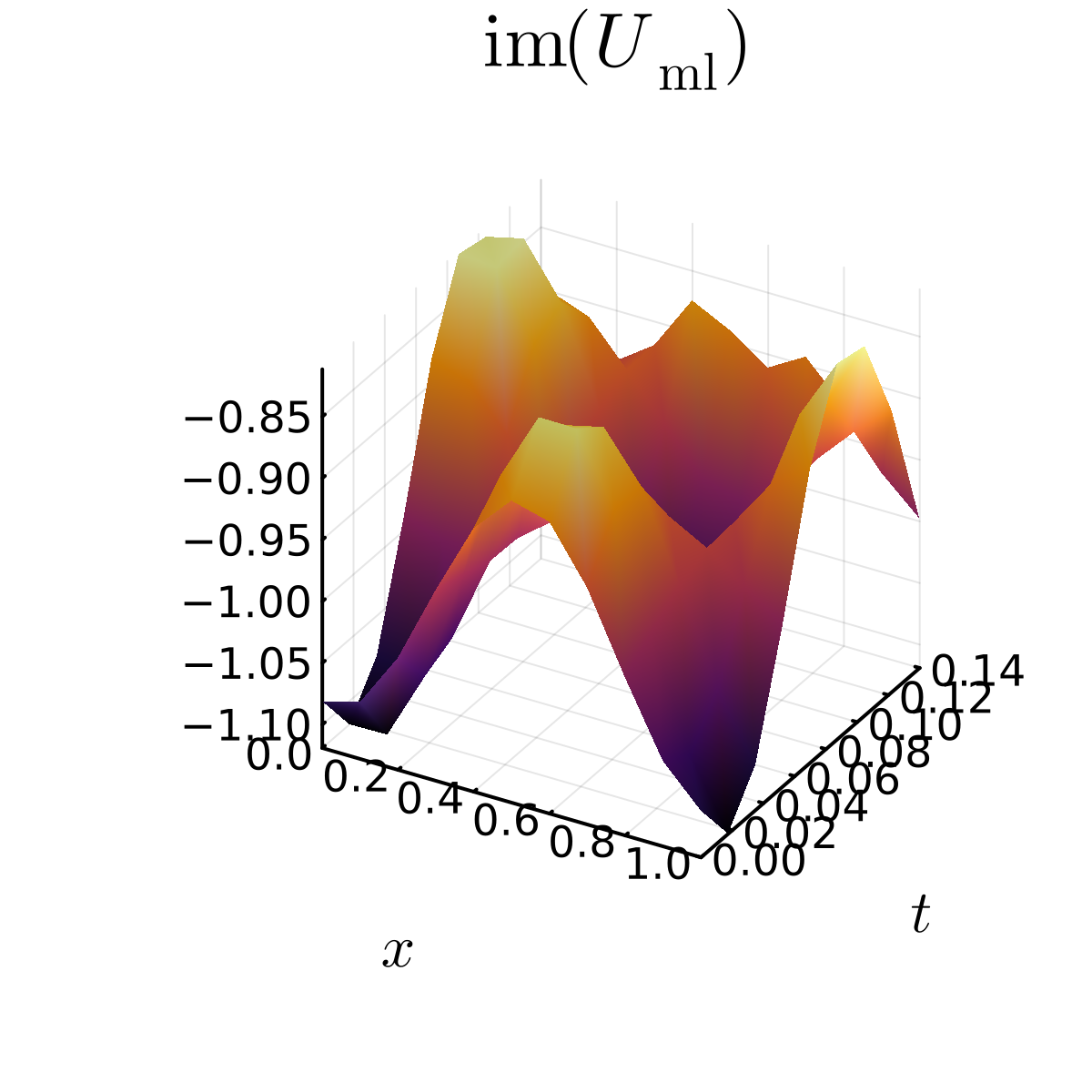}
	\includegraphics[width=0.32\linewidth]{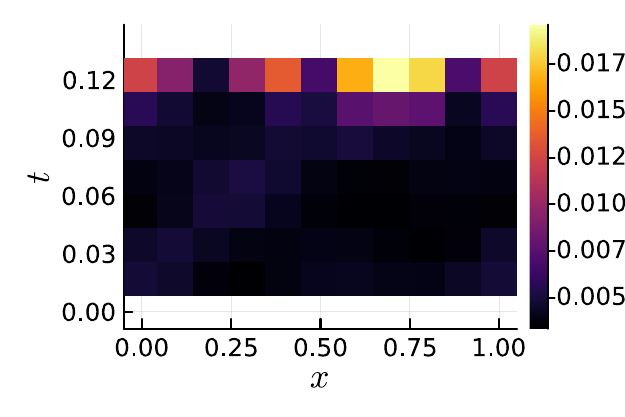}
	\caption{Schrödinger equation experiment.
		From left to right: true solution (not part of training data), predicted solution, standard deviation of $\mathrm{DEL}(\xi^M)$ at stencils of the predicted solution. First row refers to real part, second row to imaginary parts of the fields.
	}\label{fig:PredictionsSE}
\end{figure}
Moreover, \cref{fig:PredictionsSE} shows that the trained model can successfully predict motions from unseen initial data that was sampled from the same distribution as the training data.
The discrete $l_2$ norm of the error is smaller $3.9e-3$. Model uncertainty along the computed field is quantified as the component-wise standard deviation of $\mathrm{DEL}_{\mathfrak{u}^{(i,j)}}(\xi^M)$ and visualised in \cref{fig:PredictionsSE} displaying an increase in uncertainty as the solution is propagated forward in time.


\subsection*{Reproducibility}
Computations have been performed in Julia \cite{bezanson2017julia} using packages NLsolve.jl (part of the Optim.jl project \cite{mogensen2018optim}) and ForwardDiff.jl \cite{revels2016}, among others.
The source code of the experiments can be found at \url{https://github.com/Christian-Offen/Lagrangian_GP_PDE}.

\section{Convergence Analysis}\label{sec:ConvergenceAnalysis}

In this section, we will state and prove convergence results for the considered method to learn discrete Lagrangian densities from temporal-spatial data.
For notational convenience, the statement (\cref{thm:ConvergenceThmLd3Pt}) is formulated for three-point Lagrangian densities related to a 1+1-dimensional temporal-spatial dynamical systems. However, \cref{thm:ConvergenceThmLd3Pt} and its proof extend in a straight forward manner to other types of discrete Lagrangian densities such as the four-point Lagrangian densities employed in the numerical experiment of \cref{sec:ExperimentSE} or, more generally, for the class of discrete Lagrangian densities described in Section II.B of \cite{DLNNPDE}.

Additionally, \cref{app:DiscreteLagrangiansFromPDETheorem}
extracts a statement about temporal discrete Lagrangians as a special case of the presented theory and relates the statement to the discussion in our article on (temporal) Lagrangian identification \cite{DLGPode}.

\subsection{Statement of convergence theorem}

\begin{theorem}\label{thm:ConvergenceThmLd3Pt}
Let $\Omega \subset \R^d \times \R^d \times \R^d$ be open, bounded, non-empty domains. For $\mathfrak u = (u,u^+,u_+,u^-,u^-_+,u_-,u_-^+) \in (\R^d)^7$ consider the projections
$\mathrm{pr}_1(\mathfrak{u}) = (u,u^+,u_+)$,
$\mathrm{pr}_2(\mathfrak{u}) = (u^-,u,u_+^-)$,
$\mathrm{pr}_3(\mathfrak{u}) = (u_-,u_-^+ ,u)$.
Consider a sequence of observations
\[\Omega_0 = \{\mathfrak{u}^{(j)} = (u^{(j)},{u^+}^{(j)},{u_+}^{(j)},{u^-}^{(j)},{u^-_+}^{(j)},{u_-}^{(j)},{u_-^+}^{(j)}) \}_{j=1}^\infty
\subset \bigcap_{l=1}^3 \mathrm{pr}_l^{-1}(\Omega)
\]
of a dynamical system governed by the discrete Euler--Lagrange equation of an (unknown) discrete Lagrangian density $L_d^\rf \in \mathcal{C}^1(\overline{\Omega})$, i.e.\
$\mathrm{DEL}(L_d^\rf)(\mathfrak{u})=0$ for all $\mathfrak u \in \Omega_0$.
Let $K \in \mathcal{C}^1(\overline{\Omega} \times \overline{\Omega})$ be a continuously differentiable kernel on $\Omega$, $\mathcal u_b \in \Omega$, $r_b \in \R$, $p_b \in \R$ and assume that $L^\rf_d$ is contained in the reproducing kernel Hilbert space $B$ to $K$ and fulfils the normalisation condition
	\begin{equation}\label{eq:NormalisationThmLd3Pt}
		\Phi_b(L^\rf_d) = (p_b,r_b) \quad \text{with} \quad
		\Phi_b(L_d)= \left( \nabla_2 L_d(\mathcal u_b),L_d(\mathcal u_b) \right)
	\end{equation}
	and that $B$ embeds continuously into $\mathcal{C}^1(\overline{\Omega})$.
	Let $\xi \in \mathcal{N}(0,\mathcal K)$ be a canonical Gaussian random variable over $B$. Then the sequence of conditional means $L_{d,(j)}$ of $\xi$ conditioned on the first $j$ observations and the normalisation conditions
	\begin{equation}
		\mathrm{DEL}(\xi)(\mathfrak{u}^{(i)}) =0\, (\forall i\le j),\quad
		\Phi_N(\xi) = (p_b,r_b)
	\end{equation}
	converges in the reproducing kernel Hilbert space norm $\| \cdot \|_B$ and in the uniform norm of the space of differentiable functions $\| \cdot \|_{\mathcal{C}^1(\overline{\Omega})}$ to a Lagrangian $L_{d,(\infty)} \in B$ that is 
	\begin{itemize}
	\item consistent with the normalisation $\Phi_b(L_{d,(\infty)}) = (p_b,r_b)$
	\item consistent with the dynamics, i.e.\ $\mathrm{DEL}(L_{d,(\infty)})(\mathfrak{u}) =0$ for all $\mathfrak{u} \in \overline{\Omega_0}$ 
	\item and $L_{d,(\infty)}$ minimizes $\| \cdot \|_B$ among all discrete Lagrangian densities with the properties above.
\end{itemize}
\end{theorem}

\begin{Remark}
By addition of null-Lagrangian densities (see \cref{prop:normaliseLd}) $L_d^{\rf}$ can be assumed to fulfil the normalisation condition \eqref{eq:NormalisationThmLd3Pt} without loss of generality. However, it is stated in this form to assume compatibility with the reproducing kernel Hilbert space $B$.
\end{Remark}

\begin{Remark}
It makes sense to consider observations $\Omega_0$ in \cref{thm:ConvergenceThmLd3Pt} that densely fill a sufficiently large space, i.e.\ $ \Omega \subset  \mathrm{pr}_k(\overline{\Omega_0})$ for $k=1,2,3$, such that the dynamics on all of $ \Omega$ is covered.
\end{Remark}

\subsection{Setting and proof of \cref{thm:ConvergenceThmLd3Pt}}

Let $\Omega \subset \R^d \times \R^d \times \R^d$ be an open, bounded, non-empty domain.
For data $\mathfrak u = (u,u^+,u_+,u^-,u^-_+,u_-,u_-^+) \in (\R^d)^7$ define the projections
$\mathrm{pr}_1(\mathfrak{u}) = (u,u^+,u_+)$,
$\mathrm{pr}_2(\mathfrak{u}) = (u^-,u,u_+^-)$,
$\mathrm{pr}_3(\mathfrak{u}) = (u_-,u_-^+ ,u)$.
To a discrete Lagrangian $L_d^\rf \in \mathcal{C}^1(\overline\Omega)$ consider the set of admissible stencil data
\begin{equation*}\label{eq:OmegaHatLd}
\hat \Omega_{\mathrm{pre}}
= \left\{
\mathfrak{u} = (u,u^+,u_+,u^-,u^-_+,u_-,u_-^+) \, \left| \,
\mathrm{DEL}(L_d^\rf)(\mathfrak{u})=0, \mathrm{pr}_k(\mathfrak{u}) \in \overline\Omega, k=1,2,3
\right.\right\}.
\end{equation*}
The set $\hat \Omega_{\mathrm{pre}}$ is closed in $(\R^d)^7$.


\begin{remark}\label{rem:OmegaPre}
In degenerate/symmetric situations, $\hat \Omega_{\mathrm{pre}}$ can contain additional stencils that do not relate to discrete fields: for instance, for $L_d^\rf(u,u^+,u_+) = \frac {(u^+-u)^2}{2\Delta t^2} - \frac {(u_+-u)^2}{2\Delta x^2}$ the discrete Euler--Lagrange equation $\mathrm{DEL}(L_d^\rf)(\mathfrak{u})=0$ is a five-point stencil only relating $u,u^+,u^-,u_+,u_-$ rather than a seven-point stencil relating $u,u^+,u_+,u^-,u^-_+,u_-,u_-^+$. Thus
$(u,u^+,u_+,u^-,r_1,u_-,r_2) \in \Omega_{\mathrm{pre}}$ for all $r_1,r_2 \in \R$ whenever 
\[
 \frac{u^- - 2 u + u^+}{\Delta t^2}-\frac{u_- - 2 u + u_+}{\Delta x^2} =0.
\]
To cover the case where training data stencils are collected from observed discrete fields, we will, therefore, proof the statement for
observations that are dense in some closed subset $\hat \Omega$ of $\hat \Omega_{\mathrm{pre}}$.
\end{remark}

Let $\Omega_0 = \{\mathfrak{u}^{(j)}\}_{j=1}^\infty \subset \hat \Omega_{\mathrm{pre}}$ be a sequence and $\hat \Omega := \overline \Omega_0 \subset \hat \Omega_{\mathrm{pre}}$ (topological closures in $\hat \Omega_{\mathrm{pre}}$ and $(\R^d)^7$ coincide).
The discrete Lagrangian operator $\mathrm{DEL}$ constitutes a bounded linear operator
\[\Phi^{\infty} \colon \mathcal{C}^1(\overline\Omega) \to \mathcal{C}^0(\hat{\Omega},\R^d),
\quad L_d \mapsto \mathrm{DEL}(L_d).
\]

Since for each $\mathfrak{u}$ the evaluation functional $\mathrm{ev}_{\mathfrak{u}} \colon f \mapsto f(\mathfrak{u})$ on $\mathcal{C}^0(\hat \Omega,\R^d)$ is bounded, the following functions constitute bounded linear functionals for $j \in \N$:
\begin{align*}
	\Phi_j &\colon \mathcal{C}^1(\overline\Omega) \to \R^d, \quad \quad\Phi_j(L_d)=\Phi^{(\infty)} (L_d)(\mathfrak{u}^{(j)})\\
	\Phi^{(j)} &\colon \mathcal{C}^1(\overline\Omega) \to (\R^d)^j, \quad \Phi^{(j)} =(\Phi_1,\ldots,\Phi_j).
\end{align*}

For a reference point $\mathcal u_b \in \Omega$ and for $p_b \in \R^d$, $r_b \in \R$ we define the bounded linear functional
\begin{align}\label{eq:PhiNDiscrete3Pt}
	\Phi_b &\colon \mathcal{C}^1(\overline\Omega) \to \R^{d+1}, \quad \quad \Phi_b(L)= \left( \nabla_2 L_d(\mathcal u_b), L_d(\mathcal u_b) \right),
\end{align}
related to our normalisation condition on the local momentum $\mathrm{Mm}^+$ for seven-point Lagrangians (\cref{sec:Normalisation}).
Here $\nabla_2 L_d$ denotes the derivative with respect to the second component of $L_d$.
We will further use the shorthands
\[\Phi_b^{(k)} = (\Phi_1,\ldots,\Phi_k,\Phi_b) \quad \text{and}\quad \Phi_b^{(\infty)} = (\Phi^{(\infty)} ,\Phi_b),\]
and define
\begin{align*}
	y^{(k)} &= (0,\ldots,0,p_b,r_b) \in (\R^d)^k \times \R^d\times \R\\
	y^{(\infty)} &= (0,p_b,r_b) \in \mathcal{C}^0(\overline\Omega,\R^d) \times \R^d\times \R.
\end{align*}

Consider the following assumption.
\begin{Assumption}\label{ass:LinBDiscrete3Pt}
	Assume that there is a reflexive, uniformly convex Banach space $B$ with continuous embedding $B \hookrightarrow \mathcal{C}^1(\overline\Omega)$ such that
	\[\{L_d \in \mathcal{C}^1(\overline \Omega) \, | \, \Phi_b^{(\infty)}(L_d) = y^{(\infty)} \} \cap B \not = \emptyset\]
	In other words, $B$ is assumed to contain a Lagrangian consistent  with the normalisation and the underlying dynamics.
\end{Assumption}

The affine linear subspace
\begin{align*}
	A^{(j)} &= \{ L_d\in B \, | \, \Phi_b^{(j)}(L_d) = y^{(j)} \} \quad (j \in \N)\\
	A^{(\infty)} &= \{ L_d \in B \, | \, \Phi_b^{(\infty)}(L_d) = y^{(\infty)} \}
\end{align*}
are closed in $B$ and not empty by \cref{ass:LinBDiscrete3Pt}. Therefore, the following minimisations constitute convex optimisation problems on $B$ with unique minima in $A^{(j)}$ or $A^{(\infty)}$, respectively:
\begin{equation}\label{eq:MinProblemsDiscrete}
	\begin{split}
		{L_d}_{(j)} &= \arg \min_{{L_d} \in A^{(j)}} \| {L_d} \|_B \\
		{L_d}_{(\infty)} &= \arg \min_{{L_d} \in A^{(\infty)}} \| {L_d} \|_B 
	\end{split}
\end{equation}
Here $\| \cdot \|_B $ denotes the norm in $B$.

\begin{proposition}\label{prop:ConvergenceInBDiscrete3pt}
	The minima ${L_d}_{(j)}$ converge to ${L_d}_{(\infty)}$ in the norm $\| \cdot \|_B$ and, thus, in $\| \cdot \|_{\mathcal C^1(\overline{\Omega})}$.
\end{proposition}

\begin{proof}
The proof consists of the following steps and follows the same strategy as in the case of temporal Lagrangians \cite{DLGPode}. We will prove:

\begin{enumerate}
	
	\item The sequence $({L_d}_{(j)})_j$ has a weakly convergent subsequence.
	\item The weak limit $L_{d(\infty)}^\dagger$ is contained in $A^{(\infty)}$.
	\item 
	$L_{d(\infty)}^\dagger$ coincides with $L_{d(\infty)}$.
	\item The whole sequence $({L_d}_{(j)})_j$ converges weakly against $L_{d(\infty)}$.
	\item The norms $\|{L_d}_{(j)}\|_B$ converge to $\|L_{d(\infty)}\|_B$.
	\item Step 4 and 5 imply that the sequence $({L_d}_{(j)})_j$ converges strongly against $L_{d(\infty)}$.
	
\end{enumerate}



{\em Step 1.} The sequence of affine spaces $A^{(1)} \supseteq A^{(2)} \supseteq A^{(2)} \supseteq \ldots$ decreases monotonously. Moreover, $\emptyset \not = A^{(\infty)} \subseteq \bigcap_{j=1}^\infty A^{(j)}$.
It follows that the norms $\|L_{d(j)}\|_B$ increase monotonously. The sequence is bounded by $\|L_{d(\infty)}\|_B$.
By reflexivity of $B$, we find a weakly convergent subsequence $(L_{d(j_i)})_{i \in \N}$ (see \cite[Theorem 3.18]{Brezis2011}). Let $L_{j_i} \rightharpoonup L^\dagger_{d(\infty)}$ denote the weak limit in $B$.
We have
\begin{equation}\label{eq:WeakLimitUpperBoundLd3}
	\| L^\dagger_{d(\infty)} \|_B \le \liminf_{i \to \infty} \| L_{d(j_i)} \|_B \le \| L_{d(\infty)} \|_B
\end{equation}
by the weak lower semi-continuity of the norm.

{\em Step 2.} 
We show $L^\dagger_{d(\infty)} \in A^{(\infty)}$.

Let $\mathfrak u \in \hat \Omega$. As the sequence $\Omega_0 = (\mathfrak{u}^{(m)} )_{m=1}^\infty$ is dense in $\hat \Omega$, there exists a subsequence $(\mathfrak{u}^{(m_l)} )_{l=1}^\infty$ converging to $\mathfrak{u}$.

For each $l$ there is $N$ such that $j_N \ge m_l$. For $i \ge N$ it holds that $\Phi^{(\infty)}_b (L_{d(j_i)}) (\mathfrak{u}^{(m_l)})=0$.
Therefore, for each $l$
\[
\lim_{i \to \infty} \Phi^{(\infty)}_b (L_{d(j_i)}) (\mathfrak{u}^{(m_l)}) =0.\]

For each $l$ the linear operator $\Phi^{(\infty)}_b(\cdot) (\mathfrak{u}^{(m_l)}) \colon B \to \R^d \times \R^{d+1}$ into the finite-dimensional space $\R^d \times \R^{d+1}$ is bounded.
By the weak convergence $L_{d(j_i)} \rightharpoonup L_{d(\infty)}^\dagger$ we have
\begin{align*}\label{eq:UseWeakConvergenceLd3}
	\Phi^{(\infty)}_b(L^\dagger_{d(\infty)}) (\mathfrak{u}^{(m_l)})  
	= \lim_{i \to \infty} \Phi^{(\infty)}_b (L_{d(j_i)}) (\mathfrak{u}^{(m_l)})=0.
\end{align*}

By continuity $\Phi^{(\infty)}_b(L^\dagger_{d(\infty)}) \in \mathcal{C}^0(\overline{\Omega})$, we have
\begin{align*}\label{eq:UseContinuity1}
	\Phi^{(\infty)}_b(L^\dagger_{d(\infty)}) (\mathfrak{u})
	= \lim_{l \to \infty} \Phi^{(\infty)}_b(L^\dagger_{d(\infty)}) (\mathfrak{u}^{(m_l)}) =0.
\end{align*}

As this holds for all $\mathfrak{u}\in \hat \Omega$ we conclude $L^\dagger_{d(\infty)} \in A^{(\infty)}$.

{\em Step 3.} As $L^\dagger_{d(\infty)} \in A^{(\infty)}$ by Step 2, we have $\| L^\dagger_{d(\infty)} \|_B \ge \| L_{d(\infty)} \|_B$ since $L_{d(\infty)}$ is the unique minimiser of the minimisation problem of \eqref{eq:MinProblemsDiscrete}. Together with \eqref{eq:WeakLimitUpperBoundLd3} we conclude $\| L^\dagger_{d(\infty)} \|_B = \| L_{d(\infty)} \|_B$. By uniqueness of the minimiser $L_{d(\infty)}$ this implies $ L^\dagger_{d(\infty)}  =  L_{d(\infty)}$.
It follows that $L_{d(j_i)}$ converges weakly to $ L_{d(\infty)}$.

{\em Step 4.}
The choice of a weakly convergent subsequence of $L_{d(j)}$ in Step 1 was arbitrary. Thus any subsequence of $L_{d(j)}$ has a weakly convergent subsequence and any weakly convergent subsequence has weak limit $ L_{d(\infty)}$. It follows that the whole series $L_{d(j)}$ converges weakly to the unique weak limit $L_{d(\infty)}$.

{\em Step 5.}
As $L_{d(j)} \rightharpoonup L_{d(\infty)}$, and $\|L_{d(j)}\|_B \le \|L_{d(\infty)}\|_B$, and by the weak lower semi-continuity of the norm 
\begin{equation*}\label{eq:NormConvergence}
\| L_{d(\infty)} \|_B  \le \liminf_{j \to \infty} \| L_{d(j)} \|_B
\le \limsup_{j \to \infty} \| L_{d(j)} \|_B
\le \| L_{d(\infty)} \|_B.
\end{equation*}
Therefore
\begin{equation}\label{eq:NormConvergenceFinal}
	\lim_{j \to \infty} \|L_{d(j)}\|_B = \|L_{d(\infty)} \|_B.
	\end{equation}

{\em Step 6.}
As $B$ is uniformly convex, by \cite[Prop. 3.32]{Brezis2011} weak convergence 
$L_{d(j)} \rightharpoonup L^\dagger_{d(\infty)}$ together with \eqref{eq:NormConvergenceFinal} implies strong convergence $\lim_{j \to \infty} L_{d(j)} = L_{d(\infty)}$.
\end{proof}

\begin{proof}[\cref{thm:ConvergenceThmLd3Pt}]
By \cite[Thm 12.5]{OwhadiScovel2019OptimalRecoverySplines} the unique minimisers $L_{d(j)}$ in \eqref{eq:MinProblemsDiscrete} coincide with the conditional means considered in \cref{thm:ConvergenceThmLd3Pt}.
Moreover, reproducing kernel Hilbert spaces are uniformly convex, reflexive Banach spaces. Thus the statement follows directly from \cref{prop:ConvergenceInBDiscrete3pt}.
\end{proof}

\section{Summary and further directions}\label{sec:Summary}
We have introduced a method based on Gaussian Process regression to obtain data-driven models of discrete Lagrangian densities that define discrete field theories.
The data consists of discrete observations of fields of the true field theory, which we decompose into stencil data $\mathfrak U$ for the discrete Euler--Lagrange equation of a discrete Lagrangian density.
A canonical Gaussian process $\xi$ modelling a discrete Lagrangian density is conditioned to fulfil the discrete Euler--Lagrange equation  $\mathrm{DEL}_{\mathrm{u}}(\xi) =0$ on all observed stencils $\mathfrak u \in \mathfrak{U}$.
As these conditions are linear in $\xi$, the posterior process is Gaussian as well and the conditional mean and the conditional covariance operator can be explicitly computed.
To avoid degeneracies of the posterior processes, we further condition the process on non-triviality conditions. We show that this does not restrict the generality of our ansatz.
Next to providing a method to predict solutions to field theories, we use the statistical framework of Gaussian process regression to obtain a quantification of model uncertainty of the learned discrete Euler--Lagrange equations and any linear observable in the conditioned Gaussian process. 
Our method is structure-preserving in the sense that the learned discrete field theories are guaranteed to fulfil a discrete variational principle by construction.
To illustrate the method, we provide numerical examples based on the wave equation and the Schrödinger equation.

The method can be interpreted as a meshless collocation method \cite{SchabackWendland2006} for solving the discrete Euler--Lagrange equation for the discrete Lagrangian $L_d$. 
Overcoming the ambiguity of discrete Lagrangian densities (gauge freedom), we provide a statement that guarantees convergence of the posterior means to a true discrete Lagrangian density as the distance of data points converges to zero.
For this, the proof exploits a characterisation of posterior means as minimisers of a reproducing kernel Hilbert space norm constrained by the observations \cite{OwhadiScovel2019OptimalRecoverySplines}. 

In future work, it is of interest to prove convergence rates for the provided method and to develop efficient computational methods for evaluations of the machine-learned theory in high-dimensional, large data-regimes.
Moreover, techniques of this article can be combined with Lie group based methods to learn symmetric representations of (discrete) temporal Lagrangians (such as \cite{SymLNN}) which may be extended to Lagrangian densities and temporal-spatial symmetries. These are of interest for the data-driven identification of conservation laws and the detection of structurally simple solutions such as travelling \cite{DLNNDensity,DLNNPDE}.
Furthermore, rather than learning of one specific variational formulation to a field theory, learning of alternative (i.e.~non-gauge equivalent) variational formulations is of interest for system identification as alternative variational formulations allow for the derivation of a series of conservation laws {\cite{HENNEAUX198245,Marmo1987,MARMO1989389,Carinena1983}.
	Finally and more generally, it is of interest to obtain a fundamental understanding of the interplay between geometric prior knowledge, data requirements, and model uncertainty to clarify the role of geometry in machine learning.

\section*{Acknowledgments}
The author would like to thank Konstantin Sonntag for a helpful discussion that has advanced the convergence proof.
Moreover, the author acknowledges the Ministerium für Kultur und Wissenschaft des Landes Nordrhein-Westfalen and computing time provided by
the Paderborn Center for Parallel Computing (PC2).

\section*{Data availability}
The data that support the findings of this study are openly available in the GitHub repository Christian-Offen/Lagrangian\_GP\_PDE at
\url{https://github.com/Christian-Offen/Lagrangian_GP_PDE}.

\appendix

\section*{Appendices}

\section{Continuous Lagrangian theories}\label{app:ContinuousLagrangianDensities}

Let us review of continuous Euler--Lagrange dynamics. For details and a more formal treatment we refer to \cref{sec:ConvergenceAnalysis} and the literature on variational calculus \cite{gelfand2000calculus,RoubicekCalculusofVariations}.

\subsection{The Euler--Lagrange equations for continuous theories}

The Euler--Lagrange operator for a first order field theory on $\R^n$ or subsets of $\R^n$ with coordinate $t=(t_1,\ldots,t_n)$ and $\R^d$-valued field $u=(u^1,\ldots,u^n)$ is given as 
\begin{equation}
	\begin{split}
		\mathrm{EL}(L)^r&((u^s)_s,(u^s_{t_k})_{s,k},(u^s_{t_l,t_k})_{s,l,k})\\
		&= \sum_{k=1}^{n} \frac{\p }{\p t_k} \left( \frac{\p L}{\p {u_{t_k}^r}} \right) - \frac{\p L}{\p u^r}\\
		&= \sum_{k,l=1}^{n}\sum_{s=1}^d \left( \frac{\p^2 L}{\p {u_{t_k}^r}\p {u_{t_l}^s}} u^s_{t_k,t_l} \right) 
		+\sum_{k=1}^{n}\sum_{s=1}^d \left( \frac{\p^2 L}{\p {u_{t_k}^r}\p {u^s}} u^s_{t_l} \right)
		- \frac{\p L}{\p u^r}
	\end{split}
\end{equation}
for $r=1,\ldots,d$.
The Euler--Lagrange equation
\begin{equation}
	\mathrm{EL}(L)\left(u(t), \frac{\p u}{\p t_k}(t),\frac{\p u}{\p t_k \p t_l}(t)\right) = 0
\end{equation}
constitutes a system of $d$ partial differential equations.
Here we use the convention $\mathrm{EL}(L) = (\mathrm{EL}(L)^1,\ldots,\mathrm{EL}(L)^d)$.
Unless degenerate, the equations are of second order.
In case $n=1$ the equation simplifies to the following set of second order ordinary differential equation:
\begin{equation}\label{eq:ELODE}
	0=\mathrm{EL}(L)(u,u_t,u_{tt})
	= \frac{\p^2 L}{\p u \p u_t} u_t + \frac{\p^2 L}{\p u_t \p u_t} u_{tt} - \frac{\p L}{\p u}.
\end{equation}
Here the first and second summand denote matrix vector multiplications and $\frac{\p L}{\p u}$ denotes the gradient of $L$ with respect to $u$.

\subsection{Ambiguity of Lagrangian densities}

Lagrangian densities can be ambiguous in two different ways: 
\begin{enumerate}
	
	\item Lagrangians $L$ and $\tilde L$ can yield the same Euler--Lagrange operator up to rescaling, i.e.\ 
	\[
	\mathrm{EL}(L) = \rho \cdot  \mathrm{EL}(\tilde{L}), \quad \rho \in \R\setminus \{0\}
	\]
	and, therefore, the same Euler--Lagrange equations up to rescaling. We call $L$ and $\tilde L$ {\em (gauge-) equivalent}.
	For equivalent Lagrangians $L$, $\tilde L$ there exists $\rho \in \R \setminus \{0\}$, $c \in \R$ such that $L - \rho \tilde L -c$ is a total divergence
	\[
	L - \rho \tilde L -c = \mathrm{div}_t F
	\]
	for a continuously differentiable function $F=(F^1,\ldots,F^n)$, where
	\begin{equation}
		\mathrm{div}_t F(u,(u_{t_k})_k)
		= \sum_{s=1}^d \sum_{k=1}^n u_{t_k}^s \frac{\p F^j}{\p u^s} (u).
	\end{equation}
	(We have restricted ourselves to autonomous Lagrangians.) This may be seen by noting that for the corresponding action functionals $S$ and $\tilde S$ 
	the difference $S(u) - \rho S(u)$ is an integral over the boundary $\p \mathcal{M}$ if $u \colon \mathcal M \to \R^d$. However, the considered variations fix $\p \mathcal{M}$ and, thus, do not influence the stationarity of a field $u$ \cite{gelfand2000calculus}.
	
	\item Even when two Lagrangians $L$ and $\tilde L$ are not equivalent, they can yield the same set of solutions $u$, i.e.
	\[
	\forall u\colon \mathcal{M} \to \R^d : \mathrm{EL}(L)(\mathfrak{u}(t)) =0 \iff \mathrm{EL}(\tilde L)(\mathfrak{u}(t)) =0,
	\]
	where $\mathcal{M} \subseteq \R^n$ is the domain of the field theory.
	Here
	\[\mathfrak{u}(t) = \left(u(t),\left(\frac{\p u(t)}{\p t_k}\right)_k,\left(\frac{\p^2 u(t)}{\p t_k \p t_l}\right)_{k,l}\right)\]
	denotes the second jet extension of a field $u$.
	In that case $\tilde{L}$ is called an {\em alternative Lagrangian density} to $L$.
	
	\begin{example}[Affine linear motions]
		For any twice differentiable $g \colon \R^d \to \R$ with nowhere degenerate Hessian matrix $\mathrm{Hess}(g)$, the Lagrangian $L(u,u_t) = g(u_t)$ describes affine linear motions in $\R^d$:
		\[
		0 = \mathrm{EL}(L) = \mathrm{Hess}(g)(u_t)u_{tt}.
		\]
	\end{example}

\end{enumerate}

In general, the existence of alternative Lagrangian densities is related to additional geometric structure and conserved quantities of the system \cite{HENNEAUX198245,Marmo1987,MARMO1989389,Carinena1983}.
Gauge equivalence, however, is exhibited by all variational systems.

\section{Alternative convergence statement for discrete Lagrangians }\label{app:DiscreteLagrangiansFromPDETheorem}

\Cref{thm:ConvergenceThmLd3Pt} provides a blueprint for convergence results for various types of discrete Lagrangians.
In \cref{thm:ConvergenceThmLdTemporal} below we reformulate \cref{thm:ConvergenceThmLd3Pt} for discrete temporal Lagrangians. 
Discrete temporal Lagrangians have already been considered separately in \cite{DLGPode}. We will then relate the statement of \cite{DLGPode} and \cref{thm:ConvergenceThmLdTemporal}.

\begin{theorem}\label{thm:ConvergenceThmLdTemporal}
Let $\Omega \subset \R^d \times \R^d$ be open, bounded, non-empty domains.
Consider a sequence of observations
\[\Omega_0 = \{\mathfrak{u}^{(j)} = ({u^-}^{(j)},u^{(j)},{u^+}^{(j)}) \}_{j=1}^\infty\]
of a dynamical system governed by the discrete Euler--Lagrange equation of an (unknown) discrete Lagrangian density $L_d^\rf \in \mathcal{C}^1(\overline{\Omega})$, i.e.\
	\[
	\Omega_0 \subset \hat \Omega = \left\{\mathfrak{u}=(u^-,u,u^+) \, \left| \,
	\mathrm{DEL}(L_d^\rf)(\mathfrak{u})=0, \begin{array}{c}
		(u,u^+) \in \overline \Omega\\
		(u^-,u) \in \overline \Omega
	\end{array}
	\right.\right\}.
	\]
	Assume $\Omega \subset \mathrm{pr_-}\overline{\Omega_0}$ and $\Omega \subset \mathrm{pr_+}\overline{\Omega_0}$
	with $\mathrm{pr}_-(u^-,u,u^+)=(u^-,u)$ and $\mathrm{pr}_+(u^-,u,u^+)=(u,u^+)$.
	Let $K\in \mathcal{C}^1(\overline{\Omega} \times \overline{\Omega})$ be a kernel on $\Omega$, $ u_b \in \Omega$, $r_b \in \R$, $p_b \in \R$ and assume that $L^\rf_d$ is contained in the reproducing kernel Hilbert space $B$ to $K$ and fulfils the normalisation condition
	\begin{equation}\label{eq:NormalisationThmLdTemporal}
		\Phi_b(L_d^\rf) = (p_b,r_b) \quad \text{with} \quad
		\Phi_b(L_d)= \left(-\nabla_2 L_d( u_b),L_d( u_b) \right)
	\end{equation}
	and that $B$ embeds continuously into $\mathcal{C}^1(\overline{\Omega})$.
	Let $\xi \in \mathcal{N}(0,\mathcal{K})$ be a canonical Gaussian random variable over $B$. Then the sequence of posterior means $L_{d,(j)}$ of the Gaussian field $\xi$ conditioned on the first $j$ observations and the normalisation conditions
	\begin{equation}
		\mathrm{DEL}(\xi)(\mathfrak{u}^{(i)}) =0\, (\forall i\le j),\quad
		\Phi_b(\xi) = (p_b,r_b)
	\end{equation}
	converges in the reproducing kernel Hilbert space norm $\| \cdot \|_B$ and in the uniform norm of the space of differentiable functions $\| \cdot \|_{\mathcal{C}^1(\overline{\Omega})}$ to a Lagrangian $L_{d,(\infty)} \in B$ that is 
	\begin{itemize}
		\item consistent with the normalisation $\Phi_b(L_{d,(\infty)}) = (p_b,r_b)$
		\item consistent with the dynamics, i.e.\ $\mathrm{DEL}(L_{d,(\infty)})(\mathfrak{u}) =0$ for all $\mathfrak{u} \in\overline{\Omega_0}$ 
		\item and it minimizes $\| \cdot \|_B$ among discrete Lagrangians with the above properties.
	\end{itemize}
\end{theorem}

\begin{proof}
The proof of \cref{thm:ConvergenceThmLd3Pt} can be adapted easily.
\end{proof}

In comparison to \cite[Theorem 2]{DLGPode}, \cref{thm:ConvergenceThmLdTemporal} requires $\Omega=\Omega_a=\Omega_b$ (slightly more restrictive). Under the assumptions that a globally Lipschitz continuous discrete flow map $g$ exists governing the observed dynamical system,
any dense set $\{(x_0^{(j)},x_1^{(j)})\}_{j=1}^\infty \subset \Omega_a$ gives rise to a dense, countable set of observations
\[
\Omega_0 = \{(x_0^{(j)},g(x_0^{(j)},x_1^{(j)}))\} \subset \hat \Omega
\]
where $\hat \Omega = \{(x_0,g(x_0,x_1)) \, | \, (x_0,x_1) \in \overline \Omega\}$. Hence, \cite[Theorem 2]{DLGPode} follows from \cref{thm:ConvergenceThmLdTemporal} in case $\Omega_a=\Omega_b$.

\bibliographystyle{plainurl}
\bibliography{literature}

\begin{thebibliography}{10}

\bibitem{Blanchette2020}
Christine Allen-Blanchette, Sushant Veer, Anirudha Majumdar, and Naomi~Ehrich
  Leonard.
\newblock {LagNetViP}: A {L}agrangian neural network for video prediction
  ({AAAI} 2020 symposium on physics guided ai), 2020.
\newblock \href {https://doi.org/10.48550/ARXIV.2010.12932}
  {\path{doi:10.48550/ARXIV.2010.12932}}.

\bibitem{Aoshima2021}
Takehiro Aoshima, Takashi Matsubara, and Takaharu Yaguchi.
\newblock Deep discrete-time lagrangian mechanics.
\newblock {\em ICLR SimDL}, 5 2021.
\newblock URL: \url{https://simdl.github.io/files/49.pdf}.

\bibitem{bezanson2017julia}
Jeff Bezanson, Alan Edelman, Stefan Karpinski, and Viral~B Shah.
\newblock Julia: A fresh approach to numerical computing.
\newblock {\em SIAM review}, 59(1):65--98, 2017.
\newblock \href {https://doi.org/10.1137/141000671}
  {\path{doi:10.1137/141000671}}.

\bibitem{Brezis2011}
Haim Brezis.
\newblock {\em Functional Analysis, Sobolev Spaces and Partial Differential
  Equations}.
\newblock Springer New York, New York, NY, 2011.
\newblock \href {https://doi.org/10.1007/978-0-387-70914-7}
  {\path{doi:10.1007/978-0-387-70914-7}}.

\bibitem{Carinena1983}
J~F Carinena and L~A Ibort.
\newblock Non-noether constants of motion.
\newblock {\em Journal of Physics A: Mathematical and General}, 16(1):1, 1
  1983.
\newblock \href {https://doi.org/10.1088/0305-4470/16/1/010}
  {\path{doi:10.1088/0305-4470/16/1/010}}.

\bibitem{OwhadiLearningPDEGP}
Yifan Chen, Bamdad Hosseini, Houman Owhadi, and Andrew~M. Stuart.
\newblock Solving and learning nonlinear pdes with {G}aussian processes.
\newblock {\em Journal of Computational Physics}, 447:110668, 2021.
\newblock \href {https://doi.org/10.1016/j.jcp.2021.110668}
  {\path{doi:10.1016/j.jcp.2021.110668}}.

\bibitem{ChenOwhadiSChaefer2024SparseCholesky}
Yifan Chen, Houman Owhadi, and Florian Schäfer.
\newblock Sparse {C}holesky factorization for solving nonlinear pdes via
  {G}aussian processes.
\newblock {\em Mathematics of Computation}, June 2024.
\newblock \href {https://doi.org/10.1090/mcom/3992}
  {\path{doi:10.1090/mcom/3992}}.

\bibitem{ChristmannSteinwart2008RKHS}
Andreas Christmann and Ingo Steinwart.
\newblock {\em Kernels and Reproducing Kernel Hilbert Spaces}, pages 110--163.
\newblock Springer New York, New York, NY, 2008.
\newblock \href {https://doi.org/10.1007/978-0-387-77242-4_4}
  {\path{doi:10.1007/978-0-387-77242-4_4}}.

\bibitem{LNN}
Miles Cranmer, Sam Greydanus, Stephan Hoyer, Peter Battaglia, David Spergel,
  and Shirley Ho.
\newblock Lagrangian neural networks, 2020.
\newblock \href {https://doi.org/10.48550/ARXIV.2003.04630}
  {\path{doi:10.48550/ARXIV.2003.04630}}.

\bibitem{evangelisticdc2022}
Giulio Evangelisti and Sandra Hirche.
\newblock Physically consistent learning of conservative lagrangian systems
  with gaussian processes.
\newblock In {\em 2022 IEEE 61st Conference on Decision and Control (CDC)}.
  IEEE, 2022.
\newblock \href {https://doi.org/10.1109/CDC51059.2022.9993123}
  {\path{doi:10.1109/CDC51059.2022.9993123}}.

\bibitem{gelfand2000calculus}
I.M. Gelfand, S.V. Fomin, and R.A. Silverman.
\newblock {\em Calculus of Variations}.
\newblock Dover Books on Mathematics. Dover Publications, 2000.

\bibitem{HENNEAUX198245}
Marc Henneaux.
\newblock Equations of motion, commutation relations and ambiguities in the
  {L}agrangian formalism.
\newblock {\em Annals of Physics}, 140(1):45--64, 1982.
\newblock \href {https://doi.org/10.1016/0003-4916(82)90334-7}
  {\path{doi:10.1016/0003-4916(82)90334-7}}.

\bibitem{Karniadakis2021}
George~Em Karniadakis, Ioannis~G. Kevrekidis, Lu~Lu, Paris Perdikaris, Sifan
  Wang, and Liu Yang.
\newblock Physics-informed machine learning.
\newblock {\em Nature Reviews Physics}, 3(6):422--440, Jun 2021.
\newblock \href {https://doi.org/10.1038/s42254-021-00314-5}
  {\path{doi:10.1038/s42254-021-00314-5}}.

\bibitem{SymLNN}
Yana Lishkova, Paul Scherer, Steffen Ridderbusch, Mateja Jamnik, Pietro Liò,
  Sina Ober-Blöbaum, and Christian Offen.
\newblock Discrete {L}agrangian neural networks with automatic symmetry
  discovery.
\newblock {\em IFAC-PapersOnLine}, 56(2):3203--3210, 2023.
\newblock 22nd IFAC World Congress.
\newblock \href {https://doi.org/10.1016/j.ifacol.2023.10.1457}
  {\path{doi:10.1016/j.ifacol.2023.10.1457}}.

\bibitem{Lu2021}
Lu~Lu, Pengzhan Jin, Guofei Pang, Zhongqiang Zhang, and George~Em Karniadakis.
\newblock Learning nonlinear operators via {DeepONet} based on the universal
  approximation theorem of operators.
\newblock {\em Nature Machine Intelligence}, 3(3):218--229, mar 2021.
\newblock URL: \url{https://doi.org/10.1038%2Fs42256-021-00302-5}, \href
  {https://doi.org/10.1038/s42256-021-00302-5}
  {\path{doi:10.1038/s42256-021-00302-5}}.

\bibitem{MARMO1989389}
Giuseppe Marmo and G.~Morandi.
\newblock On the inverse problem with symmetries, and the appearance of
  cohomologies in classical {L}agrangian dynamics.
\newblock {\em Reports on Mathematical Physics}, 28(3):389--410, 1989.
\newblock \href {https://doi.org/10.1016/0034-4877(89)90071-2}
  {\path{doi:10.1016/0034-4877(89)90071-2}}.

\bibitem{Marmo1987}
Giuseppe Marmo and C.~Rubano.
\newblock On the uniqueness of the {L}agrangian description for charged
  particles in external magnetic field.
\newblock {\em Il Nuovo Cimento A}, 98(4):387--399, 10 1987.
\newblock \href {https://doi.org/10.1007/bf02902083}
  {\path{doi:10.1007/bf02902083}}.

\bibitem{Marsden1998}
Jerrold~E. Marsden, George~W. Patrick, and Steve Shkoller.
\newblock Multisymplectic geometry, variational integrators, and nonlinear
  pdes.
\newblock {\em Communications in Mathematical Physics}, 199(2):351–395, 12
  1998.
\newblock \href {https://doi.org/10.1007/s002200050505}
  {\path{doi:10.1007/s002200050505}}.

\bibitem{Marsden2001}
Jerrold~E. Marsden, Sergey Pekarsky, Steve Shkoller, and Matthew West.
\newblock Variational methods, multisymplectic geometry and continuum
  mechanics.
\newblock {\em Journal of Geometry and Physics}, 38(3):253--284, 2001.
\newblock \href {https://doi.org/10.1016/S0393-0440(00)00066-8}
  {\path{doi:10.1016/S0393-0440(00)00066-8}}.

\bibitem{MarsdenWestVariationalIntegrators}
Jerrold~E. Marsden and Matthew West.
\newblock Discrete mechanics and variational integrators.
\newblock {\em Acta Numerica}, 10:357--514, 2001.
\newblock \href {https://doi.org/10.1017/S096249290100006X}
  {\path{doi:10.1017/S096249290100006X}}.

\bibitem{Blanchette2022}
Justice Mason, Christine Allen-Blanchette, Nicholas Zolman, Elizabeth Davison,
  and Naomi Leonard.
\newblock Learning interpretable dynamics from images of a freely rotating 3d
  rigid body, 2022.
\newblock \href {https://doi.org/10.48550/ARXIV.2209.11355}
  {\path{doi:10.48550/ARXIV.2209.11355}}.

\bibitem{mogensen2018optim}
Patrick~Kofod Mogensen and Asbj{\o}rn~Nilsen Riseth.
\newblock Optim: A mathematical optimization package for {Julia}.
\newblock {\em Journal of Open Source Software}, 3(24):615, 2018.
\newblock \href {https://doi.org/10.21105/joss.00615}
  {\path{doi:10.21105/joss.00615}}.

\bibitem{LagrangianShadowIntegrators}
Sina Ober-Blöbaum and Christian Offen.
\newblock Variational learning of {E}uler–{L}agrange dynamics from data.
\newblock {\em Journal of Computational and Applied Mathematics}, 421:114780,
  2023.
\newblock \href {https://doi.org/10.1016/j.cam.2022.114780}
  {\path{doi:10.1016/j.cam.2022.114780}}.

\bibitem{DLGPode}
Christian Offen.
\newblock Machine learning of continuous and discrete variational odes with
  convergence guarantee and uncertainty quantification, 2024.
\newblock \href {https://arxiv.org/abs/2404.19626} {\path{arXiv:2404.19626}}.

\bibitem{DLNNDensity}
Christian Offen and Sina Ober-Bl{\"o}baum.
\newblock Learning discrete lagrangians for variational pdes from data
  and detection of travelling waves.
\newblock In Frank Nielsen and Fr{\'e}d{\'e}ric Barbaresco, editors, {\em
  Geometric Science of Information}, volume 14071, pages 569--579, Cham, 2023.
  Springer Nature Switzerland.
\newblock \href {https://doi.org/10.1007/978-3-031-38271-0_57}
  {\path{doi:10.1007/978-3-031-38271-0_57}}.

\bibitem{DLNNPDE}
Christian Offen and Sina Ober-Blöbaum.
\newblock Learning of discrete models of variational {PDE}s from data.
\newblock {\em Chaos}, 34:013104, 1 2024.
\newblock \href {https://doi.org/10.1063/5.0172287}
  {\path{doi:10.1063/5.0172287}}.

\bibitem{OwhadiScovel2019}
Houman Owhadi and Clint Scovel.
\newblock {\em Operator-Adapted Wavelets, Fast Solvers, and Numerical
  Homogenization: From a Game Theoretic Approach to Numerical Approximation and
  Algorithm Design}.
\newblock Cambridge Monographs on Applied and Computational Mathematics.
  Cambridge University Press, 2019.
\newblock \href {https://doi.org/10.1017/9781108594967}
  {\path{doi:10.1017/9781108594967}}.

\bibitem{OwhadiScovel2019OptimalRecoverySplines}
Houman Owhadi and Clint Scovel.
\newblock {\em Optimal Recovery Splines}, page 154–159.
\newblock Cambridge Monographs on Applied and Computational Mathematics.
  Cambridge University Press, 2019.
\newblock \href {https://doi.org/10.1017/9781108594967.017}
  {\path{doi:10.1017/9781108594967.017}}.

\bibitem{Qin2020}
Hong Qin.
\newblock Machine learning and serving of discrete field theories.
\newblock {\em Scientific Reports}, 10(1), 11 2020.
\newblock \href {https://doi.org/10.1038/s41598-020-76301-0}
  {\path{doi:10.1038/s41598-020-76301-0}}.

\bibitem{QuinRasmussen2005}
Joaquin Qui{{\~n}}onero-Candela and Carl~Edward Rasmussen.
\newblock A unifying view of sparse approximate gaussian process regression.
\newblock {\em Journal of Machine Learning Research}, 6(65):1939--1959, 2005.
\newblock URL: \url{http://jmlr.org/papers/v6/quinonero-candela05a.html}.

\bibitem{revels2016}
Jarrett Revels, Miles Lubin, and Theodore Papamarkou.
\newblock Forward-mode automatic differentiation in {J}ulia, 2016.
\newblock \href {https://arxiv.org/abs/1607.07892} {\path{arXiv:1607.07892}}.

\bibitem{RoubicekCalculusofVariations}
Tomáš Roubíček.
\newblock {\em Calculus of Variations}, pages 1--38.
\newblock John Wiley \& Sons, Ltd, 2015.
\newblock \href {https://doi.org/10.1002/3527600434.eap735}
  {\path{doi:10.1002/3527600434.eap735}}.

\bibitem{Rudy2017}
Samuel~H. Rudy, Steven~L. Brunton, Joshua~L. Proctor, and J.~Nathan Kutz.
\newblock Data-driven discovery of partial differential equations.
\newblock {\em Science Advances}, 3(4):e1602614, 2017.
\newblock \href {https://doi.org/10.1126/sciadv.1602614}
  {\path{doi:10.1126/sciadv.1602614}}.

\bibitem{SchabackWendland2006}
Robert Schaback and Holger Wendland.
\newblock Kernel techniques: From machine learning to meshless methods.
\newblock {\em Acta Numerica}, 15:543–639, 2006.
\newblock \href {https://doi.org/10.1017/S0962492906270016}
  {\path{doi:10.1017/S0962492906270016}}.

\bibitem{Schaeffer2017}
Hayden Schaeffer.
\newblock Learning partial differential equations via data discovery and sparse
  optimization.
\newblock {\em Proceedings of the Royal Society A: Mathematical, Physical and
  Engineering Sciences}, 473(2197):20160446, 2017.
\newblock \href {https://doi.org/10.1098/rspa.2016.0446}
  {\path{doi:10.1098/rspa.2016.0446}}.

\bibitem{Schaefer2021}
Florian Sch\"{a}fer, Matthias Katzfuss, and Houman Owhadi.
\newblock Sparse {C}holesky factorization by {K}ullback--{L}eibler
  minimization.
\newblock {\em SIAM Journal on Scientific Computing}, 43(3):A2019--A2046, 2021.
\newblock \href {https://doi.org/10.1137/20M1336254}
  {\path{doi:10.1137/20M1336254}}.

\bibitem{Sharma2022LagrangianROM}
Harsh Sharma and Boris Kramer.
\newblock Preserving {L}agrangian structure in data-driven reduced-order
  modeling of large-scale dynamical systems, 2022.
\newblock \href {https://doi.org/10.48550/ARXIV.2203.06361}
  {\path{doi:10.48550/ARXIV.2203.06361}}.

\bibitem{sharma2023symplectic}
Harsh Sharma, Hongliang Mu, Patrick Buchfink, Rudy Geelen, Silke Glas, and
  Boris Kramer.
\newblock Symplectic model reduction of {H}amiltonian systems using data-driven
  quadratic manifolds, 2023.
\newblock \href {https://arxiv.org/abs/2305.15490} {\path{arXiv:2305.15490}}.

\bibitem{Sharma2022}
Harsh Sharma, Zhu Wang, and Boris Kramer.
\newblock Hamiltonian operator inference: Physics-preserving learning of
  reduced-order models for canonical {H}amiltonian systems.
\newblock {\em Physica D: Nonlinear Phenomena}, 431:133122, 2022.
\newblock \href {https://doi.org/10.1016/j.physd.2021.133122}
  {\path{doi:10.1016/j.physd.2021.133122}}.

\end{thebibliography}

\end{document}